\documentclass[11pt]{amsart}
\usepackage{amssymb}

\usepackage{fixltx2e}
\usepackage[utf8]{inputenc}
\usepackage{amssymb, latexsym,amsthm,amsfonts, amsbsy, amsmath, mathrsfs}
\usepackage{bm}                          
\usepackage{verbatim}                    
\usepackage[all, knot]{xy}               
        \xyoption{arc}
        \xyoption{web}                   
\makeatletter                            
\def\MT@register@subst@font{\MT@exp@one@n\MT@in@clist\font@name\MT@font@list
 \ifMT@inlist@\else\xdef\MT@font@list{\MT@font@list\font@name,}\fi}
\makeatother

\newtheorem{thm}{Theorem}[section]

\newtheorem{lem}[thm]{Lemma}
\newtheorem{cor}[thm]{Corollary}
\newtheorem{remk}[thm]{Remark}

\newtheorem{defn}[thm]{Definition}

\newtheorem{exmp}[thm]{Example}

\newcommand{\al}{\alpha}

\newcommand{\rig}{\rightarrow}
\newcommand{\mrig}{\mathrel{-\!\!\!\!\!\rightarrow}}

\newcommand{\Rig}{\Rightarrow}

\newcommand{\bcdw}{\mathbin{\boldsymbol\cdot}}

\newcommand{\seteq}{\mathrel{\mbox{\,:\!}=\nolinebreak }\,}

\newcommand{\simop}{\mathop{{\sim\!}}}
\newcommand{\up}{\mathop{\,\uparrow}}
\newcommand{\down}{\mathop{\,\downarrow}}
\newcommand{\bup}{\mathop{\,\boldsymbol{\uparrow}}}

\newcommand{\sbA}{{\boldsymbol{A}}}
\newcommand{\sbB}{{\boldsymbol{B}}}
\newcommand{\sbC}{{\boldsymbol{C}}}
\newcommand{\sbD}{{\boldsymbol{D}}}
\newcommand{\sbE}{{\boldsymbol{E}}}

\newcommand{\sbG}{{\boldsymbol{G}}}
\newcommand{\sbH}{{\boldsymbol{H}}}

\newcommand{\sbP}{{\boldsymbol{P}}}
\newcommand{\sbQ}{{\boldsymbol{Q}}}

\newcommand{\sbW}{{\boldsymbol{W}}}
\newcommand{\sbX}{{\boldsymbol{X}}}
\newcommand{\sbY}{{\boldsymbol{Y}}}

\newcommand{\sbe}{{\boldsymbol{e}}}
\newcommand{\sbf}{{\boldsymbol{f}}}

\newcommand{\sbx}{{\boldsymbol{x}}}
\newcommand{\sby}{{\boldsymbol{y}}}
\newcommand{\sbz}{{\boldsymbol{z}}}

\let\oper=\mathbb                                  
\bmdefine{\A}{A}                                   
\bmdefine{\B}{B}
\bmdefine{\C}{C}
\bmdefine{\D}{D}
\bmdefine{\Fm}{Fm}                                 
\bmdefine{\Con}{\textup{Con}}                      

\newcommand{\VVV}{\oper{V}}

\bmdefine{\boldstar}{\mathchoice{\textstyle*}{\textstyle*}{\textstyle*}{\scriptstyle*}}

\bmdefine{\leibniz}{\varOmega}                       

\bmdefine{\btau}{\tau}                               
\bmdefine{\brho}{\rho}                               

\theoremstyle{theorem}
\newtheorem{Theorem}{Theorem}[section]
\newtheorem{Lemma}[Theorem]{Lemma}
\newtheorem{Corollary}[Theorem]{Corollary}
\newtheorem{Claim}[]{Claim}

\theoremstyle{definition}

\begin{document}

\title[Epimorphisms
in Varieties of Residuated Structures]{Epimorphisms
in Varieties of\\ Residuated Structures}
\author[G.\ Bezhanishvili]{Guram Bezhanishvili}
\address{Department of Mathematical Sciences, New Mexico State University, Las Cruces NM 88003, USA}
\email{guram@math.nmsu.edu}
\author[T.\ Moraschini]{Tommaso Moraschini}
\address{Institute of Computer Science, Academy of Sciences of the Czech Republic, Pod Vod\'{a}renskou v\v{e}\v{z}\'{i} 2, 182 07 Prague 8, Czech Republic}
\email{moraschini@cs.cas.cz}
\author[J.G.\ Raftery]{James Raftery}
\address{Department of Mathematics and Applied Mathematics,
 University of Pretoria,
 Private Bag X20, Hatfield,
 Pretoria 0028, South Africa}
\email{james.raftery@up.ac.za}
\keywords{Epimorphism, Brouwerian algebra, Heyting algebra, Esakia space, residuated lattice, Sugihara monoid,
substructural logic, intuitionistic logic, relevance logic, R-mingle, Beth definability.\vspace{1mm}
\vspace{1mm}}
\thanks{The second author acknowledges grant GA17-04630S of the Czech Grant Agency.
\\
\indent
The third author was supported in part by the National Research Foundation of South Africa (UID 85407).}

\begin{abstract}
It is proved that epimorphisms are surjective in a range of varieties of residuated structures, including all varieties of
Heyting or Brouwerian algebras of finite depth, and all varieties consisting of G\"{o}del algebras, relative Stone algebras,
Sugihara monoids or positive Sugihara monoids.
This establishes the infinite deductive Beth definability property for a
corresponding range of substructural logics.  On the other hand, it is shown that epimorphisms need not be surjective in a
locally finite variety of Heyting or Brouwerian algebras of width $2$.
It follows that the infinite Beth property is strictly stronger than the
so-called finite Beth property, confirming a conjecture of Blok and Hoogland.
\end{abstract}

\maketitle

\makeatletter
\renewcommand{\labelenumi}{\text{(\theenumi)}}
\renewcommand{\theenumi}{\roman{enumi}}
\renewcommand{\theenumii}{\roman{enumii}}
\renewcommand{\labelenumii}{\text{(\theenumii)}}
\renewcommand{\p@enumii}{\theenumi(\theenumii)}
\makeatother

{\allowdisplaybreaks

\section{Introduction}\label{introduction}

A morphism $h$ in a category $\mathsf{C}$ is called a ($\mathsf{C}$--)\,{\em epimorphism\/} provided that,
for any two $\mathsf{C}$--morphisms $f,g$ from the co-domain of $h$ to a single object,
\[
\textup{if $f\circ h=g\circ h$, then $f=g$.}
\]
We shall not distinguish notationally between a class $\mathsf{K}$ of similar algebras and the concrete category
of algebraic homomorphisms between its members.  Clearly, in such a category, every surjective $\mathsf{K}$--morphism
is a $\mathsf{K}$--epimorphism.  If the converse holds, then $\mathsf{K}$ is said to have the {\em epimorphism
surjectivity property}, or briefly, the {\em ES property}.

This property fails, for instance, in the variety of rings.  There,
the inclusion $\mathbb{Z}\mrig\mathbb{Q}$ is an epimorphism, mainly because multiplicative
inverses are `implicitly defined',
i.e., uniquely determined {\em or\/} non-existent.  The failure of surjectivity reflects the absence of an
explicit unary term defining inversehood in the language of rings.
In a slogan: epimorphisms correspond to implicit definitions and
surjective homomorphisms to explicit ones.

Groups, modules over a given ring, semilattices
and lattices
each form a variety in which
all epimorphisms are surjective; see the references in \cite{KMPT83}.
The ES property need not persist in subvarieties, however.  Indeed, it fails for {\em distributive\/}
lattices, where an embedding of the three-element chain in a four-element Boolean lattice is an
epimorphism (owing to the uniqueness of existent complements).

As this suggests, it is generally difficult to determine whether epimorphisms are surjective in a given variety.
Here, for a range of varieties of residuated structures, we shall prove that they are.
The ES property is algebraically natural, but our main motivation comes from logic, as residuated structures algebraize
substructural logics \cite{GJKO07}.

The algebraic counterpart $\mathsf{K}$ of an algebraizable logic $\,\vdash$
is a {\em pre\-variety}, i.e., a class of similar algebras, closed under isomorphisms, subalgebras and direct products;
see \cite{BH06,BP89,Cze01,CP99,FJP03}.
In this situation,
\begin{quote}
$\mathsf{K}$
has the ES property iff $\,\vdash$ has the {\em infinite (deductive) Beth
(definability) property\/} \cite[Thm.\,3.17]{BH06}.
\end{quote}
The latter signifies that, in $\,\vdash$, whenever a set $Z$ of variables is defined
{\em implicitly\/}
in terms of a disjoint set $X$ of variables by means of some set $\Gamma$ of formulas
over $X\cup Z$, then $\Gamma$ also defines $Z$ {\em explicitly\/} in terms of $X$.
In substructural logics,
this means, more precisely, that whenever
\begin{equation}\label{implicit}
\Gamma\cup \sigma
[\Gamma]\vdash z\leftrightarrow \sigma(z)
\end{equation}
holds for all $z\in Z$ and all uniform substitutions $\sigma$ (of formulas for variables) satisfying $\sigma(x)=x$ for all $x\in X$,
then
for each $z\in Z$, there is a formula $\varphi_z$ over $X$ only, such that
\begin{equation}\label{explicit}
\Gamma\vdash z\leftrightarrow\varphi_z.
\end{equation}
Here, $X,Z$ and $\Gamma$ may be infinite; no bound on their cardinalities is assumed.  Formulas in the range of
$\sigma|_{X\cup Z}$ may also involve arbitrarily many variables beyond $X\cup Z$.
To make sense of (\ref{explicit}), we assume that $X\neq\emptyset$, unless there are constant symbols in the signature.

The {\em finite Beth property\/} makes the same demand, but only when
$Z$ is finite---or equivalently,
as it turns out, when $Z$ is a singleton \cite[Cor.\,3.15]{BH06}.
\begin{exmp}\label{cpl}
\textup{(\cite{BH06})\,
In classical propositional logic ($\mathbf{CPL}$), and in
its implication fragment ($\mathbf{CPL}_{\,\rig}$), if $X=\{x_1,x_2\}$ and $Z=\{z\}$ and
\[
\Gamma=\{z\rig x_1,\,\,z\rig x_2,\,\,x_1\rig(x_2\rig z)\},
\]
then (\ref{implicit}) holds for every substitution $\sigma$
fixing $x_1$ and $x_2$.  In $\mathbf{CPL}$, (\ref{explicit}) is witnessed as
$\Gamma\vdash_\mathbf{CPL} z\leftrightarrow (x_1\wedge x_2)$,
but there is demonstrably no such instantiation in $\mathbf{CPL}_{\,\rig}$.
(Equivalently, the algebras for $\mathbf{CPL}_{\,\rig}$ need not be meet semilattice-ordered, but existent greatest lower bounds are unique.)
This shows that
$\mathbf{CPL}_{\,\rig}$
lacks even the finite Beth property, whereas $\mathbf{CPL}$ has the infinite Beth property, because epimorphisms are surjective in its algebraic counterpart---the variety of
Boolean algebras.\footnote{\,For the historical origins of the strong amalgamation (and hence the ES) property in Boolean algebras, see \cite[Footnote~7, p.\,336]{Pig72}.
\,E.W.\ Beth's original definability theorems for
classical propositional and predicate logic
were proved in \cite{Bet53}.}}\qed
\end{exmp}
\noindent
Strictly speaking, it is $\,\vdash_\mathbf{CPL}$ that is algebraized by Boolean algebras and that has the infinite Beth property, but we routinely attribute to a formal
system $\mathbf{F}$ the significant properties of its deducibility relation $\,\vdash_\mathbf{F}$.

An algebraizable logic has the finite Beth property iff its algebraic counterpart has
the `weak' ES property defined below.
(Again, see \cite{BH06}; a restricted form of this claim, due to I.\ N\'{e}meti, appeared earlier in \cite[Thm.\,5.6.10]{HMT7185}.)
It is pointed out in \cite{BH06} that the meaning of the weak ES property would not change if we allowed finite sets to play the role of the singleton
$\{b\}$ in Definition~\ref{weak es def}.
\begin{defn}\label{weak es def}
\textup{A homomorphism $h\colon\sbA\mrig\sbB$ between algebras is {\em almost-onto\/} if $\sbB$ is generated by $h[A]\cup\{b\}$ for some $b\in B$.
A prevariety $\mathsf{K}$ has the {\em weak ES property\/} if every almost-onto $\mathsf{K}$--epimorphism is surjective.}
\end{defn}

In \cite{BH06}, the Beth properties are formulated more generally---for logics that are `equivalential' in the sense of \cite{Cze01}.
Even in that wide context, it was not previously established whether the finite Beth property implies the infinite one.  A negative answer was
conjectured by Blok and Hoogland in \cite[p.\,76]{BH06}.

We shall confirm their conjecture here, by exhibiting a variety with the weak ES
property but not the ES property, algebraizing a fairly orthodox logic.
(Actually, any prevariety with the weak ES but not the ES property would confirm
the conjecture, as these properties are categorical in prevarieties and every prevariety is categorically equivalent
to one that algebraizes a sentential logic \cite[Thm.\,6.26]{Mor}.)

Rings and distributive lattices do not assist us here, as they lack even the weak ES property.  So do modular lattices, by
\cite[Thm.\,3.3]{Fre79} and its proof.
In seeking the counter-example, we must avoid
amalgamable prevarieties, because of the
following result, which combines observations in \cite{Isb66,KMPT83,Rin72} and \cite[Sec.\,2.5.3]{Hoo01}.  (Definitions of the pertinent amalgamation
properties can be found, for instance, in \cite[p.\,3204]{GR15}; they will not be needed here.)
\begin{thm}\label{hoogland et al}
A prevariety\/ $\mathsf{K}$ has the
amalgamation and weak ES properties
iff it has the strong amalgamation property.

In that case, it has
the following\/
\textup{`strong ES property':}
whenever\/ $\sbA$ is a subalgebra of some\/
$\sbB\in\mathsf{K}$ and\/ $b\in B\smallsetminus A$\textup{,} then there are two homomorphisms from\/ $\sbB$ to a single member of\/ $\mathsf{K}$ that
agree on\/ $A$ but not at\/ $b$\textup{.}
\end{thm}

Looking to models of intuitionistic logic, we recall that
the weak ES property holds in every variety of Heyting or Brouwerian algebras \cite{Kre60}.
In both cases, there are uncountably many such varieties, and Maksimova \cite{Mak00,Mak03} has shown that only finitely many of them have
the strong ES property.  It is therefore sensible
to ask which varieties of Heyting or Brouwerian algebras
have surjective epimorphisms.

Using Esakia duality, we prove that every variety of Heyting or Brouwerian algebras {\em of finite depth\/} has the ES property (Theorems~\ref{new main}
and \ref{main brouwerian}).  At depth $3$, this already supplies
$2^{\aleph_0}$ examples where the ES property holds but the strong one fails.  Another consequence is that epimorphisms are surjective in every
{\em finitely generated\/} variety of Heyting or Brouwerian algebras.
Exploiting category equivalences in \cite{GR12,GR15}, we then obtain the ES
property for a range of
varieties of non-integral residuated structures, including
all varieties consisting of Sugihara monoids or positive Sugihara monoids; see Sections~\ref{non-integral varieties section} and \ref{sugihara monoids section}.
(The results accommodate the models of various relevance logics and/or many-valued logics.)

Nevertheless, we show that epimorphisms need not be surjective in a {\em locally finite\/} variety of Heyting or Brouwerian algebras (Theorem~\ref{Thm : ES fails}).
This affirms that the
infinite Beth property is strictly stronger than the finite one, even for locally tabular logics of a long-established
kind.

For additional information about definability in substructural (and other) logics, the reader may consult
\cite{GM05,Hoo00,Hoo01,KO10}.

\section{Residuated Structures}\label{residuated structures section}

An algebra $\sbA=\langle A;\bcdw,\rig,\wedge,\vee,\sbe\rangle$ is called a {\em commutative residuated lattice},
or briefly a {\em CRL}, if $\langle A;\wedge,\vee\rangle$ is a lattice
and $\langle A;\bcdw,\sbe\rangle$ is a commutative monoid, while $\rig$ is a binary operation such that $\sbA$ satisfies
\[
x\bcdw y\leqslant z\;\Longleftrightarrow\;x\leqslant y\rig z,
\]
where $\leqslant$ is the lattice order (cf.\ \cite{GJKO07}).
We call
$\sbA$ {\em idempotent\/} if $a\bcdw a=a$ for all $a\in A$, {\em distributive\/} if its lattice reduct is distributive, and
{\em integral\/} if $\sbe$ is its greatest element.

A {\em bounded CRL\/} is the expansion of a CRL by a distinguished element $\bot$, which is the least element of the order, whence $\top\seteq\bot\rig\bot$
is the greatest element.  In an integral bounded CRL, therefore, $\sbe=\top$.  Even in an integral unbounded CRL, we tend to write $\sbe$ as $\top$, and
we have $a\leqslant b$ iff $a\rig b=\top$ (whereas in an arbitrary CRL, $a\leqslant b$ iff $\sbe\leqslant a\rig b$).

A {\em deductive filter\/} of a (possibly bounded) CRL $\sbA$ is a lattice filter of $\langle A;\wedge,\vee\rangle$ that is also a submonoid
of $\langle A;\bcdw,\sbe\rangle$.  The lattice of deductive filters of $\sbA$ and the congruence lattice ${\boldsymbol{\mathit{Con}}}\,\sbA$ of $\sbA$
are isomorphic.  The isomorphism and its inverse are given by
\begin{eqnarray*}
& F\,\mapsto\,\leibniz F\seteq\{\langle a,b\rangle\in A^2\colon a\rig b,\,b\rig a\in F\};\\
& \theta\,\mapsto\,\{a\in A\colon \langle a\wedge \sbe, \sbe\rangle\in\theta\}.
\end{eqnarray*}
We abbreviate $\sbA/\leibniz F$ as $\sbA/F$.  It follows that $\sbA$ is finitely subdirectly irreducible (i.e., its identity relation is meet-irreducible
in ${\boldsymbol{\mathit{Con}}}\,\sbA$) iff its smallest deductive filter $\{a\in A\colon \sbe\leqslant a\}$ is meet-irreducible in its lattice of deductive
filters, and that amounts to the join-irreducibility of $\sbe$ in $\langle A;\wedge,\vee\rangle$ (a condition called `well-connectedness' in much of the literature).  If $\textup{$\{a\in A\colon a<\sbe\}$}$ has a greatest
element, then $\sbA$ is subdirectly irreducible; the converse holds when $\sbA$ is idempotent.  See, for instance, \cite[Cor.\,14]{GOR08} and \cite[Thm.\,2.4]{OR07}.

If a CRL $\sbA$ is both integral and idempotent, then its operations $\bcdw$ and $\wedge$ coincide
and $\langle A;\rig,\wedge,\vee,\top\rangle$ is called a {\em Brouwerian algebra}, in which case it is distributive and determined by its lattice reduct.
In these algebras, deductive filters are just lattice filters (the latter are
assumed non-empty here);
they will be referred to simply as `filters'.

A {\em Heyting algebra\/} is a bounded Brouwerian algebra.  Thus, $\bot$ belongs to its subalgebras, and homomorphisms between Heyting algebras preserve $\bot$.
In Heyting algebras, $\neg a$ abbreviates $a\rig\bot$.

CRLs form a variety that algebraizes a rich fragment of linear logic.
The varieties of Heyting algebras algebraize the {\em super-intuitionistic logics},
i.e., the axiomatic extensions of the intuitionistic propositional logic $\mathbf{IPL}$.
The axiomatic extensions of $\mathbf{IPL}$'s negation-less fragment (a.k.a.\ the {\em positive super-intuitionistic logics\/})
are algebraized by the varieties of Brouwerian algebras.  All of these varieties are congruence distributive, as their members
have lattice reducts.

In the next theorem, the first assertion follows from a logical argument of Kreisel \cite{Kre60}, in view of the correspondences discussed
in the introduction (or see \cite[Sec.\,12]{GR15}).  For Heyting algebras, the second was proved in \cite{EG81}; see \cite{Mak00,Mak03} for more comprehensive results.

\begin{thm}\label{heyting es}\
\begin{enumerate}
\item\label{heyting es1}
Every variety consisting of Brouwerian or Heyting algebras has the weak ES property.
\item\label{heyting es2}
The variety of all Brouwerian algebras and the variety of all Heyting algebras each have the strong ES property.
\end{enumerate}
\end{thm}

\section{Esakia Duality}

Our analysis of epimorphisms in varieties of Brouwerian or Heyting algebras will exploit Esakia duality \cite{Esa74}, so we recall some prerequisites here.

In a partially ordered set $\langle X;\leqslant\rangle$,
we define
$\up{x}=\{y\in X\colon x\leqslant y\}$
and $\up U=\bigcup_{u\in U}\up{u}$, for $U\cup\{x\}\subseteq X$, and if $U=\up U$, we call $U$ an {\em up-set\/} of $\langle X;\leqslant\rangle$.
We define $\down x$ and $\down U$
dually.

An {\em Esakia space\/} $\sbX=\langle X;\tau,\leqslant\rangle$ comprises a partially ordered set $\langle X;\leqslant\rangle$ and
a Stone space $\langle X;\tau\rangle$ (i.e., a compact Hausdorff space in which each open set is a union of clopen sets), such that
\begin{enumerate}
\item
$\up x$ is closed, for all $x\in X$, and

\item $\down U$ is clopen, for every clopen $U\subseteq X$.
\end{enumerate}
In this case, the {\em Priestley separation axiom\/} of \cite{Pri70} holds: for any $x,y\in X$,
\[
\textup{if $x\not\leqslant y$, then there is a clopen up-set $U\subseteq X$ with $x\in U$ and $y\notin U$.}
\]
Esakia spaces form a category $\mathsf{ESP}$ in which the morphisms from $\sbX$ to $\sbY$ are the so-called {\em Esakia morphisms}, i.e., the isotone continuous functions $f\colon X\mrig Y$
with the following property:
\begin{equation}\label{esakia morphism}
\textup{if $x\in X$ and $f(x)\leqslant y\in Y$, then $y=f(z)$ for some $z\in\up x$.}
\end{equation}
In other words, the Esakia morphisms $\sbX\mrig\sbY$ are the continuous functions $f$ such that $\up{f(x)}=f[\up{x}]$
(alternatively, such that $\down{f^{-1}[\{x\}]}=f^{-1}[\down{x}]$) for all $x\in X$.
In $\mathsf{ESP}$, isomorphisms are just bijective Esakia morphisms, because any continuous bijection from a compact
topological space to a Hausdorff space has a continuous inverse.

We denote by $\mathsf{HA}$ the variety (and the category) of all Heyting algebras.
The next result was established by Esakia \cite[Thm.\,3, p.\,149]{Esa74}.
\begin{thm}\label{esakia}
The categories\/ $\mathsf{HA}$ and\/ $\mathsf{ESP}$ are dually equivalent, i.e., there is a category equivalence between\/ $\mathsf{HA}$ and the opposite
category of\/ $\mathsf{ESP}$\textup{.}
\end{thm}

The contravariant functor $(-)_*
\colon\mathsf{HA}\mrig\mathsf{ESP}$ works as follows.  For each Heyting algebra $\sbA$, let
$\textup{Pr}\,\sbA$ denote the set of all prime proper filters of $\sbA$.
(A filter of $\sbA$ is {\em prime\/} if its complement is closed under the binary operation $\vee$.  To
unify our account of duality for Heyting and Brouwerian algebras, we are adopting the unusual convention
that the improper filter $A$ is prime, but in the Heyting case, $A\notin\textup{Pr}\,\sbA$.)
For each $a\in A$, let $\varphi(a)$ denote
$\{F\in\textup{Pr}\,\sbA\colon a\in F\}$ and $\varphi(a)^c$ its complement $\{F\in\textup{Pr}\,\sbA\colon a\notin F\}$.
The {\em dual space\/} $\sbA_*$
of $\sbA$ is the Esakia space $\langle \textup{Pr}\,\sbA;\tau,\subseteq\rangle$, where $\tau$ is the topology on $\textup{Pr}\,\sbA$
with subbasis
$\{\varphi(a)\colon a\in A\}\,\cup\,\{\varphi(a)^c\colon a\in A\}$.
Given a homomorphism $f\colon\sbA\mrig\sbB$ between Heyting algebras, its $(-)_*
$--image $f_*
\colon\sbB_*
\mrig
\sbA_*
$ is the $\mathsf{ESP}$--morphism $F\mapsto f^{-1}[F]$ ($F\in\textup{Pr}\,\sbB$).

In the other direction, the contravariant functor $(-)^*
\colon\mathsf{ESP}\mrig\mathsf{HA}$ sends an Esakia space $\sbX=\langle X;\tau,\leqslant\rangle$ to
$\sbX^*
\seteq\langle \textup{Cu}\,\sbX;\rig,\cap,\cup,X,\emptyset\rangle\in\mathsf{HA}$,
where $\textup{Cu}\,\sbX$ is the set of clopen up-sets of $\sbX$, on which
$U\rig V\seteq X\smallsetminus\down{(U\smallsetminus V)}$.
If $g\colon\sbX\mrig \sbY$ is a morphism in $\mathsf{ESP}$, then $g^*
\colon\sbY^*
\mrig\sbX^*$ is the Heyting algebra
homomorphism $U\mapsto g^{-1}[U]$ ($U\in\textup{Cu}\,\sbY$).

For $\sbA\in\mathsf{HA}$ and $\sbX\in\mathsf{ESP}$, the canonical isomorphisms $\sbA\cong{\sbA_*}^*$ and $\sbX\cong{\sbX^*}_*$
are given by $a\mapsto\varphi(a)$ and $x\mapsto\{U\in\textup{Cu}\,\sbX:x\in U\}$, respectively.

Given an Esakia space $\langle X;\tau,\leqslant\rangle$, if $\langle X;\leqslant\rangle$ has a greatest element $m$, then
the expansion $\sbX=\langle X;\tau,\leqslant,m\rangle$ will be called a {\em pointed Esakia space}.
In this case, we use $\textup{Cu}\,\sbX$ to denote the set of all
{\em non-empty\/} clopen up-sets of $\sbX$, i.e., all clopen up-sets to which $m$ belongs.  Let $\mathsf{PESP}$ be the category of pointed Esakia
spaces, where the morphisms between objects are just the Esakia morphisms between their unpointed reducts.  Note that any such morphism preserves
the distinguished element, in view of (\ref{esakia morphism}).

Let $\mathsf{BRA}$ denote the variety of all Brouwerian algebras.  For
$\sbA\in\mathsf{BRA}$, we now use $\textup{Pr}\,\sbA$ to stand for the set of all prime filters of $\sbA$, {\em including\/}
the improper filter $A$.  There are two ways to prove the following result.  One is an
adaptation of the
proof of
Theorem~\ref{esakia}; the other is explained below.
\begin{thm}\label{pointed esakia}
The categories\/ $\mathsf{BRA}$ and\/ $\mathsf{PESP}$ are dually equivalent.
\end{thm}
The
contravariant
functors establishing this theorem are again denoted by $(-)_*$ and $(-)^*$.
Formally, they (and the canonical isomorphisms) are defined
just
as in the case of $\mathsf{HA}$ and $\mathsf{ESP}$, but with the new definitions of $\textup{Pr}\,\sbA$ and $\textup{Cu}\,\sbX$ in place of
the old ones.  Of course, $A$ becomes the distinguished element of $\sbA_*$,
while
$\sbX^*$
no longer has a distinguished
least element (and sometimes has no least element).

For a Brouwerian algebra $\sbA$, we denote by $\sbA_\bot$ the unique Heyting algebra whose lattice reduct is got by adding a new least
element $\bot$ to $\langle A;\wedge,\vee\rangle$.  Note that
$A$ is a prime filter of $\sbA_\bot$, so the
dual of $\sbA_\bot$ is the Esakia space reduct of the dual of $\sbA$.
\begin{remk}\label{bot congruences}
\emph{For $\sbA,\sbB\in\mathsf{BRA}$, no homomorphism from $\sbA_\bot$ into a nontrivial Heyting algebra can send
an element $a\in A$ to $\bot$ (otherwise, its kernel would identify $\neg a=\bot$ with $\neg\bot=\top$).
Thus, the restrictions to $\sbA$ of the $\mathsf{HA}$--morphisms $\sbA_\bot\mrig\sbB_\bot$ are just
the $\mathsf{BRA}$--morphisms $\sbA\mrig\sbB$.}\qed
\end{remk}
There is therefore a category isomorphism from $\mathsf{BRA}$ to the full subcategory $\mathsf{C}$ of
$\mathsf{HA}$ consisting of algebras of the form $\sbA_\bot$.  And the functors $(-)_*$ and $(-)^*$
between $\mathsf{HA}$ and $\mathsf{ESP}$ restrict to a dual category equivalence between $\mathsf{C}$ and the
full subcategory of $\mathsf{ESP}$ comprising the (reducts of) pointed Esakia spaces.  In this way,
Theorem~\ref{pointed esakia} can be seen as a corollary, rather than an analogue, of Theorem~\ref{esakia}.

Note that, in both theorems, an object is finite iff its dual is.

An {\em Esakia subspace\/} (briefly, an {\em E-subspace\/}) of an Esakia space $\sbX$ is a closed up-set of $\sbX$, equipped with the restricted order
and the subspace topology.  Equivalently, it is an Esakia space whose elements belong to $X$, where the
inclusion map is both a topological embedding and an Esakia morphism \cite[Ch.\,III, Lem.\,4.11]{Esa85}.
(The terms `generated subframe' and `generated subspace' are common synonyms.)
Thus, the restriction of an Esakia morphism to an E-subspace is still an Esakia morphism.
The {\em E-subspaces\/} of an object in
$\mathsf{PESP}$ are the non-empty E-subspaces of its unpointed reduct.
E-subspaces correspond dually to homomorphic images in a sense made precise in Lemma~\ref{easy}.

The {\em disjoint union\/} $\sbX$ of finitely many Esakia spaces $\sbX_1,\dots,\sbX_n$ is their order-disjoint and topologically disjoint union,
so $X=\bigcup_{i=1}^n(X_i\times\{i\})$ and a subset $Y$ of $X$ is open iff $\{x\in X_i\colon (x,i)\in Y\}$ is open in
$\sbX_i$ for each $i$.
This is an Esakia space.

Recall that a morphism $h$ in a category $\mathsf{C}$ is called a $\textup{($\mathsf{C}$--)}$\,{\em monomorphism\/} provided that, for
any $\mathsf{C}$--morphisms $f,g$ from a single object to the domain of $h$, if $h\circ f=h\circ g$, then $f=g$.  Injective
morphisms in a concrete category are clearly monomorphisms.  The converse holds in any prevariety (as these
include their $1$-generated free algebras).  Obviously, a dual category equivalence sends
epimorphisms to monomorphisms and vice versa.

For each subvariety $\mathsf{K}$ of $\mathsf{BRA}$ or $\mathsf{HA}$,
let $\mathsf{K}_*$ denote the category of $\mathsf{K}$--{\em spaces\/}---i.e., isomorphic images of duals of algebras in
$\mathsf{K}$---equipped with all Esakia morphisms between these objects.

The next lemma is essentially contained in \cite{Esa85}.  Brief explanatory comments have been appended, because of the
limited accessibility of \cite{Esa85}.

\begin{lem}\label{easy}\
\begin{enumerate}
\item\label{Lemma : Fin Subdir Irr}
A Brouwerian or nontrivial Heyting algebra
is finitely subdirectly irreducible iff its dual has a least element\/ \textup{(}viz.\ $\{\top\}$\textup{).}

\smallskip

\item\label{easy 1}
A homomorphism\/ $h$ between Brouwerian or Heyting algebras
is surjective iff\/ $h_*$ is injective.  Also, $h$ is injective iff\/ $h_*$ is surjective.

\smallskip

\item\label{easy 1.5}
The image of a morphism in\/ $\mathsf{ESP}$ or\/ $\mathsf{PESP}$ is
an E-subspace of the co-domain.

\smallskip

\item\label{upspace lem}
The dual of an E-subspace of a (possibly pointed) Esakia space\/ $\sbX$ is a homomorphic image of\/ $\sbX^*$\textup{.}

\smallskip

\item\label{Lemma : Direct Images}
If\/ $h\colon\sbA\mrig\sbB$ is a homomorphism between Brouwerian or Heyting algebras, then there is an Esakia space isomorphism from\/
$h[\sbA]_*$ onto\/ $h_*[\sbB_*]$\textup{,} defined by\/ $F\mapsto h^{-1}[F]$\textup{.}

\smallskip

\item\label{disjoint union}
The disjoint union of finitely many Esakia spaces\/ $\sbX_1,\dots,\sbX_n$
is isomorphic to the dual of
the direct product of the Heyting algebras $\sbX^*_1,\dots,\sbX^*_n$\textup{.}

\smallskip

\item\label{easy 2}
A variety\/ $\mathsf{K}$ of Brouwerian or Heyting algebras has the ES property iff monomorphisms in\/ $\mathsf{K}_*$ are injective.
\end{enumerate}
\end{lem}
\noindent
Item~(\ref{Lemma : Fin Subdir Irr}) follows from what was said about well-connectedness in Section~\ref{residuated structures section}
(or see \cite[Thm.\,2.9]{BB08}).
In the first assertion of (\ref{easy 1}), the forward implication is easy.
Conversely, if $h_*$ is injective, then it's a monomorphism in $\mathsf{PESP}$
or $\mathsf{ESP}$, hence $h$ is an epimorphism in $\mathsf{BRA}$ or $\mathsf{HA}$, and is thus onto, by
Theorem~\ref{heyting es}(\ref{heyting es2}).
In the second assertion of (\ref{easy 1}), the forward implication instantiates the {\em prime filter
extension theorem\/}:
if $\sbA$ is a sublattice of a distributive lattice $\sbB$, then the prime filters of $\sbA$ are just its
intersections with the prime filters of $\sbB$.
The reverse implication is easy.
Item~(\ref{easy 1.5}) follows from (\ref{esakia morphism}), because a continuous function from a compact
space to a Hausdorff space sends closed sets to closed sets.
Items (\ref{upspace lem}) and (\ref{easy 2}) follow from (\ref{easy 1}).  In (\ref{easy 2}), the restrictions of
$(-)_*$ and $(-)^*$ to $\mathsf{K}$ and $\mathsf{K}_*$
establish a duality between $\mathsf{K}$ and $\mathsf{K}_*$.
In (\ref{Lemma : Direct Images}), if $\textup{$g\colon\sbA\mrig h[\sbA]$}$ is the surjective homomorphism got by
restricting the range of $h$, then $g_*\colon h[\sbA]_*\mrig\sbA_*$ is an injective Esakia morphism, by (\ref{easy 1}), but its range is $h_*[\sbB_*]$, again
by the prime filter extension theorem.  The proof of (\ref{disjoint union})
is as for Boolean algebras or bounded distributive lattices, cf.\ \cite[Lem.\,IV.4.8]{BS81}.

\section{Depth}

Let $\sbP=\langle P;\leqslant\rangle$ be a poset, with $p\in P$ and
$0<n\in\omega$.  We say that $p$ has {\em depth\/ $n$ in\/ $\sbP$} if there is a
chain $p=p_1<p_2< \ldots <p_n$ in $\sbP$ and there is no chain $p=q_1<\ldots<q_n<q_{n+1}$ in $\sbP$.
If there is no positive integer $n$ for which this holds, then
$p$ is said to have {\em depth\/} $\infty$ {\em in\/} $\sbP$.

If $p$ is an element of an Esakia space $\sbX=\langle X;\tau,\leqslant\rangle$, we define the {\em depth
of\/ $p$ in\/ $\sbX$}
to be the depth of $p$ in $\langle X;\leqslant\rangle$.

In a pointed Esakia space $\sbX=\langle X;\tau,\leqslant,m\rangle$, if $m\neq p\in X$,
it is convenient to define the
{\em depth
of\/ $p$ in\/ $\sbX$} as
the depth of $p$ in the poset $\langle X\smallsetminus\{m\};\leqslant\rangle$, declaring the {\em depth of\/} $m$ {\em in\/} $\sbX$ to be $0$.
Thus, in the passage to the Esakia space reduct of a pointed space, the finite depths of elements rise by one.

The {\em depth\/} (a.k.a.\ `height') of a non-empty Esakia space or a pointed one
is defined as the supremum of the depths of its elements.
The empty Esakia space is said to have depth $0$.  (Our assignments of depth differ by one from certain analogous traditions, such as
Krull dimension in rings, but they are convenient for
results like Theorems~\ref{heyting depth formulas} and \ref{brouwerian depth formulas} below.)

By (\ref{esakia morphism}), if the domain of a morphism in $\mathsf{ESP}$ or $\mathsf{PESP}$ has depth at most $n\in\omega$,
then so has the image.

The {\em depth\/} of a Heyting or Brouwerian algebra is defined as the depth of its dual space.
Note that a Heyting algebra and its Brouwerian reduct have the same depth.

For $n\in\omega$, we denote by $\mathsf{BRA}_n$ and by $\mathsf{HA}_n$
the respective classes of all Brouwerian and of all Heyting algebras that have depth at most $n$.
The following result was essentially proved by Ono \cite{Ono70} and by Maksimova \cite{Mak72}.
\begin{thm}\label{heyting depth formulas}
For each\/
$n\in\omega$\textup{,} a Heyting algebra
has depth at most\/ $n$ iff it satisfies\/
$h_n\thickapprox\top$\textup{,} where\/ \textup{$h_0\seteq y$} and, for\/ $n>0$\textup{,} \textup{
\begin{eqnarray*}
&& h_{n}\seteq x_n\vee(x_n\rig h_{n-1}).
\end{eqnarray*}}\noindent Consequently, $\mathsf{HA}_n$ is
a variety.\footnote{\,Ono used
implicational terms
$p_n\seteq ((x_n\rig p_{n-1})\rig x_n)\rig x_n$ instead of $h_n$.  Mak\-simova used $bd_0\seteq\bot$ and $bd_n\seteq x_n\vee(x_n\rig bd_{n-1})$.
By induction, a Heyting algebra $\sbA$ satisfies $bd_n\leqslant h_n$, but $bd_n$ instantiates $h_n$, so
$\sbA$ satisfies
$bd_n\thickapprox\top$ iff it satisfies $h_n\thickapprox\top$.  The connection between depth and the `slices' of Hosoi \cite{Hos67} was
established in \cite{Ono70}.}
\end{thm}

For a Brouwerian algebra $\sbA$, recall that $(\sbA_\bot)_*$
is the Esakia space reduct of $\sbA_*$.
Thus, $\sbA$ has finite depth iff $\sbA_\bot$ does, in which case the depth of $\sbA$ is one
less than that of $\sbA_\bot$.
\begin{thm}\label{bot lem}
For each\/ $n\in\omega$\textup{,} a Brouwerian algebra\/ $\sbA$ satisfies\/ $h_n\thickapprox\top$ iff\/
$\sbA_\bot$ satisfies\/ $h_{n+1}\thickapprox\top$\textup{.}
\end{thm}
\begin{proof}
($\Leftarrow$) \,Let $\varepsilon$ be the instance of $h_{n+1}\thickapprox\top$ in which $y$ is replaced by $\bot$.
If $a\in A$, then $a\vee (a\rig\bot)=a\vee\bot=a$.  Substituting $x_1$ for $x_1\vee(x_1\rig\bot)$ in $\varepsilon$, we get the
instance of $h_n\thickapprox\top$ in which $y,x_1,\dots,x_n$ are replaced, respectively, by $x_1,\dots,x_{n+1}$.  Thus, $\sbA$ satisfies $h_n\thickapprox\top$.

($\Rig$)
\,If any one of $x_1,\dots,x_{n+1}$ is evaluated as $\bot$ in $\sbA_\bot$,
then $h_{n+1}$ takes the value $\top$, regardless of the values in $\sbA_\bot$ of the other variables.  So, consider a valuation
$f$ in $\sbA_\bot$ for which $x_1,\dots,x_{n+1}$ receive values in $\sbA$.  If $y$ receives a value in $\sbA$, then $h_n$ takes
the value $\top$, by assumption, whence so does $h_{n+1}$.  On the other hand, if $y$
gets the value $\bot$, then $x_1\vee(x_1\rig y)$ gets the same value as $x_1$.  But, formally replacing $x_1\vee(x_1\rig y)$ by
$x_1$ in $h_{n+1}$, we obtain an instance of $h_n$ whose variables are given values in $\sbA$ by $f$, so the corresponding value
of $h_{n+1}$ is $\top$, by assumption.
\end{proof}

\begin{thm}\label{brouwerian depth formulas}
For each\/
$n\in\omega$\textup{,} a Brouwerian algebra\/ $\sbA$ has depth at most\/ $n$ iff it satisfies\/
$h_n\thickapprox\top$\textup{.}
Thus, $\mathsf{BRA}_n$
is a variety.
\end{thm}
\begin{proof}
Let $\sbX=\sbA_*$.  Then $\sbA$ and $\sbX$
have depth at most $n$ iff $\sbA_\bot$ and the Esakia space reduct of $\sbX$
have depth at most $n+1$,
iff $\sbA_\bot$ satisfies $h_{n+1}\thickapprox\top$ (Theorem~\ref{heyting depth formulas}), iff $\sbA$ satisfies $h_n\thickapprox\top$ (Theorem~\ref{bot lem}).
\end{proof}

The {\em depth\/}
of a variety $\mathsf{K}$ of Brouwerian or Heyting algebras is the supremum of the depths of its members.
If that supremum is finite, then
$\mathsf{K}$ is said to have {\em finite depth}.
If all members of $\mathsf{K}$ have finite depth, then so does $\mathsf{K}$.
This follows from Theorems~\ref{heyting depth formulas} and \ref{brouwerian depth formulas} and \L os' Theorem \cite[Thm.\,V.2.9]{BS81},
because varieties are closed under ultraproducts.

Nontrivial Boolean algebras have depth $1$.  Every finitely generated variety of Brouwerian or Heyting algebras has finite depth.  This fact will be needed
in Sections~\ref{ES in HA and BRA section}, \ref{non-integral varieties section} and \ref{sugihara monoids section},
so we give a uniform proof below, although the Heyting case is well known (see the remarks before Lemma~\ref{another bot lem}).

Here and subsequently, $\mathbb{I}$, $\mathbb{H}$, $\mathbb{S}$, $\mathbb{P}$ and $\mathbb{P}_\mathbb{U}$ stand
for closure under isomorphic and homomorphic images, subalgebras, direct products and ultraproducts, respectively,
while $\mathbb{V}$ denotes varietal closure, i.e., $\mathbb{V}=\mathbb{HSP}$.
Recall that $\mathbb{P}_\mathbb{U}(\mathsf{K})\subseteq\mathbb{I}(\mathsf{K})$
for any finite set $\mathsf{K}$ of finite algebras.
Given a class $\mathsf{M}$ of algebras, we denote by $\mathsf{M}_\textup{FSI}$
the class of all finitely subdirectly irreducible members of $\mathsf{M}$.
\,{\em J\'{o}nsson's Theorem\/} \cite{Jon67,Jon95} asserts that
if $\mathsf{M}$ is contained in a congruence distributive variety, then
$\mathbb{V}(\mathsf{M})_\textup{FSI}\subseteq \mathbb{HSP}_\mathbb{U}(\mathsf{M})$.

\begin{lem}\label{finite}
Let\/ $\mathsf{K}$ be a finitely generated variety of Brouwerian or Heyting algebras.  Then there is an integer\/ $n$
such that, whenever $F$ is a prime filter of some\/ $\sbA\in\mathsf{K}$\textup{,}
then at most\/ $n$ filters of\/ $\sbA$ contain\/ $F$\textup{.}  In particular, $\mathsf{K}$ has depth at most\/ $n$\textup{.}
\end{lem}
\begin{proof}
Let $\sbB$ be a finite
algebra with $\mathsf{K}=\mathbb{V}(\sbB)$.
Let $n$ be the
number of equivalence relations on $B$.  Let $F\in\textup{Pr}\,\sbA$, where $\sbA\in\mathsf{K}$.
Then $\leibniz F$ is meet-irreducible in ${\boldsymbol{\mathit{Con}}}\,\sbA$, i.e., $\sbA/F\in\mathsf{K}_\textup{FSI}$.
By J\'{o}nsson's Theorem
and since $\sbB$ is finite, $\sbA/F\in\mathbb{HSP}_\mathbb{U}(\sbB)=\mathbb{HS}(\sbB)$, so $\left|A/F\right|\leq\left|B\right|$.  Consequently,
$\left|{{\mathit{Con}}}(\sbA/F)\right|\leq n$, and
${\boldsymbol{\mathit{Con}}}(\sbA/F)$
is isomorphic to the interval $[\leibniz F,A^2]$ of ${\boldsymbol{\mathit{Con}}}\,\sbA$, by the Correspondence Theorem.  Applying $\leibniz^{-1}$, we see that the
interval $[F,A]$ in the filter lattice of $\sbA$ also has at most $n$ elements.
\end{proof}

\begin{lem}\label{can add bounds}
Let\/ $\mathsf{K}$ be a quasivariety of Brouwerian algebras.  If\/ $\sbA\in\mathsf{K}$\textup{,} then\/ $\sbA\in\mathbb{S}(\sbB)$ for some\/ $\sbB\in\mathsf{K}$
such that\/ $\sbB$ has a least element.
\end{lem}
\begin{proof}
Let $\mathsf{K}_\mathsf{b}$ comprise
the algebras in $\mathsf{K}$ that have a least element.  For each $\sbA\in\mathsf{K}$ and $a\in A$, the set $\up{a}$
is the universe of a subalgebra of $\sbA$ belonging to $\mathsf{K}_\mathsf{b}$, and $\sbA$ is the directed union of these subalgebras.
Thus, $\mathsf{K}$ is contained in the quasi\-variety $\mathbb{ISPP}_\mathbb{U}(\mathsf{K}_\mathsf{b})$
generated by $\mathsf{K}_\mathsf{b}$, because quasivarieties are closed under directed unions (as they are axiomatizable
by quasi-identities of finite length).  But $\mathsf{K}_\mathsf{b}$ is clearly closed under $\mathbb{P}_\mathbb{U}$, $\mathbb{P}$ and
$\mathbb{I}$, so $\mathsf{K}=\mathbb{S}(\mathsf{K}_\mathsf{b})$.
\end{proof}

It is known that each of the varieties $\mathsf{HA}_n$ is locally finite \cite{Kom75,Kuz74}.  From this we can infer:
\begin{thm}
For each\/
$n\in\omega$\textup{,} the variety\/ $\mathsf{BRA}_n$ is locally finite, i.e., every finitely generated Brouwerian algebra of finite depth is finite.
\end{thm}
\begin{proof}
Let $\sbA\in\mathsf{BRA}_n$.
Lemma~\ref{can add bounds}
supplies
a Heyting algebra $\sbB^+$, with Brouwerian reduct
$\sbB
\in\mathsf{BRA}_n
$, where $\sbA\in\mathbb{S}(\sbB
)$.
Then $\sbB^+\in\mathsf{HA}_n$.  If $\sbA$ is generated by a finite set $Y$, then $Y$ and $A$ generate the same subalgebra $\sbA'$ of $\sbB^+$,
and $\sbA'\in\mathsf{HA}_n$, so it is finite, whence $\sbA$ is finite.
\end{proof}

A locally finite subvariety of $\mathsf{HA}$ or $\mathsf{BRA}$ need not have finite depth (as will become evident in
Section~\ref{a counter-example section}).  On the other hand,
a variety $\mathsf{K}$ of Heyting algebras is finitely generated iff it has finite depth and {\em finite width}, where the latter
means that there is a finite bound on the cardinalities of antichains in the dual of any member of $\mathsf{K}_\textup{FSI}$, cf.\
\cite[Thm.\,5.1]{BGMM06}.  As with depth, the bound can be chosen uniform.
The same statements for subvarieties of $\mathsf{BRA}$ can be inferred---for instance, by applying Lemma~\ref{can add bounds} and the next result.

\begin{lem}\label{another bot lem}
Let\/ $\mathsf{K}$ be a variety of Brouwerian algebras,
and\/
$\mathsf{C}$ the class of all Heyting algebras of the form $\sbA_\bot$ such that $\sbA\in\mathsf{K}$\textup{.}
Then the nontrivial algebras in\/ $\mathbb{V}(\mathsf{C})_\textup{FSI}$ belong to\/ $\mathsf{C}$\textup{.}
\end{lem}
\begin{proof}
Let $\Sigma$ be an equational base for $\mathsf{K}$.
Then $\mathsf{C}$ is the elementary class of nontrivial Heyting algebras $\sbB$ such that $B\smallsetminus\{\bot\}$
is closed under $\wedge$ and $\rig$ and $\langle B\smallsetminus\{\bot\};\rig,\wedge,\vee,\top\rangle$ satisfies $\Sigma$.  Thus, $\mathsf{C}$ is closed under $\mathbb{P}_\mathbb{U}$,
and obviously also under $\mathbb{S}$, while nontrivial members of $\mathbb{H}(\mathsf{C})$ lie in $\mathsf{C}$, by
Remark~\ref{bot congruences}.  The lemma is therefore a consequence of J\'{o}nsson's Theorem.
\end{proof}

\section{Epimorphisms in Heyting and Brouwerian Varieties}\label{ES in HA and BRA section}

Given an Esakia space $\sbX$,
let $\textup{max}\,\sbX$ comprise the elements of depth $1$ in $\sbX$.
It is proved in \cite[Ch.\,III, Thms.\,2.3 \& 2.1]{Esa85} that
\begin{equation}\label{max x}
\textup{$\textup{max}\,\sbX$ is a closed subset of $\sbX$,}
\end{equation}
and if $x\in X$, then $x\leqslant y$ for some $y\in\textup{max}\,\sbX$.
In fact, if $\sbX=\sbA_*$, where $\sbA\in\mathsf{HA}$ and $D=\{a\in A:\neg a=\bot\}$, then $\textup{max}\,\sbX=\bigcap_{a\in D}\varphi(a)$.
When $\sbX$ has finite depth, (\ref{max x}) can be strengthened:
\begin{lem}\label{guram}
\textup{(\cite[Lem.\,7]{Bez00})}\, Let\/ $\sbX$ be an Esakia space
of finite depth.  Then, for each positive integer\/ $n$\textup{,}
the set $P$
of all elements of depth less than\/ $n$ in\/ $\sbX$ is closed in\/
$\sbX$\textup{,}
whence it is the universe of an E-subspace\/ $\sbP$ of\/ $\sbX$\textup{.}
\end{lem}

\begin{lem}\label{claim b}
Let\/ $h\colon\sbX\mrig \sbY$ be an\/ $\mathsf{ESP}$--morphism.  Let\/ $2\leq n\in\omega$ and let\/
$P$ comprise the elements of depth less than\/ $n$ in\/ $\sbX$\textup{.}
Suppose\/ $h|_P$ is injective.  If\/ $a\in X$ has depth\/ $n$ and\/ $b\in P$ and\/ $h(a)=h(b)$\textup{,} then\/ $a<b$\textup{.}
\end{lem}
\begin{proof}
Under the given assumptions, $a$ has a cover $c\in P$, because $n\geq 2$.  Now $h(b)=h(a)\leqslant h(c)$,
so $h(c)=h(d)$ for some $d\in\up{b}$, by (\ref{esakia morphism}).  Then $d\in P$.  Suppose $a\not< b$.  Then $c\neq b$, and so $c\neq d$ (as $c$ and $d$ have different depths, unless $d=b$).
By the injectivity of $h|_P$, therefore, $h(c)\neq h(d)$, a contradiction.
Thus, $a<b$, as claimed.
\end{proof}

We can now prove our first main result about the ES property.
\begin{thm}\label{new main}
Let\/ $\mathsf{K}$ be a variety of Heyting algebras, where\/ $\mathsf{K}$ has finite depth.  Then epimorphisms in\/ $\mathsf{K}$ are
surjective.
\end{thm}
\begin{proof}
The proof is by induction on the depth, $n$ say, of $\mathsf{K}$.  The result is trivial when $n=0$, i.e., when $\mathsf{K}$ is a trivial
variety.
Let $n>0$ and assume that every subvariety of $\mathsf{HA}_{n-1}$ has the ES property.
Equivalently, by Lemma~\ref{easy}(\ref{easy 2}), monomorphisms are injective in the category $\mathsf{M}_*$ of $\mathsf{M}$--spaces,
for every subvariety $\mathsf{M}$ of $\mathsf{HA}_{n-1}$.

It is likewise enough to show that monomorphisms are injective in $\mathsf{K}_*$.  And for this, it suffices to prove the following
more specialized claim.
\begin{Claim}\label{claim 0}
\,If\/ $h\colon\sbX\mrig\sbY$ is a
monomorphism in\/ $\mathsf{K}_*$\textup{,} with\/ $X=\up{\{x,y\}}$ for some\/ $x,y\in X$ such that\/
$h(x)=h(y)$\textup{,} then $x=y$\textup{.}
\end{Claim}
\noindent
To see
that this suffices, suppose $h\colon\sbX\rig\sbY$ is an arbitrary monomorphism in $\mathsf{K}_*$ (so $\sbX^*,\sbY^*\in\mathsf{K}$),
with $x,y\in X$.
In $\sbX$, the up-set $Q\seteq\up{\{x,y\}}$ is closed,
so it is the universe of an E-subspace $\sbQ$ of $\sbX$, hence the inclusion map $j\colon\sbQ\mrig\sbX$ is an Esakia morphism.
By Lemma~\ref{easy}(\ref{upspace lem}),
$\sbQ^*\in\mathbb{H}(\sbX^*)\subseteq\mathsf{K}$, so $h|_Q=h\circ j\colon\sbQ\mrig\sbY$ is
a $\mathsf{K}_*$--morphism.
As $h$ and $j$ are both $\mathsf{K}_*$--monomorphisms, so is their composition $h|_Q$ (to which Claim~\ref{claim 0} applies).

It therefore remains only to prove Claim~\ref{claim 0}.
Let $h\colon\sbX\mrig \sbY$ be a mono\-morphism
in $\mathsf{K}_*$, where $X=\up{\{x,y\}}$
and $h(x)=h(y)$.

By Lemma~\ref{guram}, $\sbX$ has an E-subspace $\sbP$, comprising
the elements of depth less than $n$ in $\sbX$.
By Lemma~\ref{easy}(\ref{upspace lem}),
$\sbP^*\in\mathbb{H}(\sbX^*)\subseteq
\mathsf{K}\cap\mathsf{HA}_{n-1}$.
By Lemma~\ref{easy}(\ref{easy 1.5}), $h[P]$ is the universe of an E-subspace $h[\sbP]$ of $\sbY$, which has depth less than $n$, since $\sbP$ has.
And $h[\sbP]^*\in\mathbb{H}(\sbY^*)\subseteq
\mathsf{K}\cap\mathsf{HA}_{n-1}$, by Lemma~\ref{easy}(\ref{upspace lem}).   In $(\mathsf{K}\cap\mathsf{HA}_{n-1})_*$, the map $h|_P\colon\sbP\mrig h[\sbP]$ is
a monomorphism, as $h$ is a monomorphism in $\mathsf{K}_*$.
Applying the induction hypothesis to $\mathsf{K}\cap\mathsf{HA}_{n-1}$, we see that
$h|_P$ is injective.

Let $\bup{\sbx}$ and $\bup{\sby}$ denote the respective E-subspaces of $\sbX$ on $\up{x}$ and $\up{y}$.

Suppose, with a view to contradiction, that $x\neq y$.  We shall contradict the fact that $h$ is a $\mathsf{K}_*$--monomorphism by constructing a space $\sbW\in\mathsf{K}_*$ and
distinct $\mathsf{K}_*$--morphisms $g_1,g_2\colon\sbW\mrig\sbX$ such that $h\circ g_1=h\circ g_2$.

As $h|_P$ is injective, one of $x,y$ has depth $n$.  This produces two cases.

\smallskip

{\bf Case~I.}  \,$x$ and $y$ both have depth $n$.

\smallskip

In this case, $x$ and $y$ have the same covers in $\sbX$.  For suppose, on the contrary, that
$y<u\in P$, where $x\not\leqslant u$.  Then $n\geq 2$ and,
as $h$ is an Esakia morphism with $h(x)=h(y)\leqslant h(u)$, we have $h(u)=h(w)$ for some $w\in\up{x}$, so $u\neq w$.  As $h|_P$ is injective, we cannot have $w\in P$, so $x=w$.
Then $h(x)=h(u)$, but this
contradicts Lemma~\ref{claim b}, because $x\not< u$.  By symmetry, therefore, $x$ and $y$ have the same covers.

Thus, $\{x\}$ and $\{y\}$ are clopen in $\bup{\sbx}$ and $\bup{\sby}$, respectively, by the Priestley separation axiom (or since $\sbX$ is Hausdorff and $(\up{x})\smallsetminus\{x\}=(\up{y})\smallsetminus\{y\}=P$,
which is closed in $\sbX$).  It follows that $\bup{\sbx}$ and $\bup{\sby}$ are isomorphic Esakia spaces.

Let $\bup{\sbz}$ (with universe $\up{z}$) be an isomorphic copy
of $\bup{\sbx}$, disjoint from $X$.  In the disjoint union $\sbW$ of $\bup{\sbx}$, $\bup{\sby}$ and $\bup{\sbz}$, we may identify the minimal elements with $x,y,z$.
Each strict upper bound $a$ of $x$ in $\sbX$ gives rise
to three copies of itself in $\sbW$, say $a_x>x$, $a_y>y$ and $a_z>z$.
Let $g_1\colon\sbW\mrig\sbX$ be the function sending $a_x,a_y,a_z$ back to $a$ whenever $x<a\in X$, where $g_1(y)=y$ and $g_1(x)=x=g_1(z)$.
It is easily checked that $g_1$ is an Esakia morphism.
By symmetry, so is the function $g_2\colon\sbW\mrig\sbX$ which differs from $g_1$ only in that $g_2(z)=y$.
Figure~\ref{fig:1} summarizes the situation.
Now $g_1\neq g_2$, but $h\circ g_1=h\circ g_2$, because $h(x)=h(y)$.

{\small
\begin{figure}[htbp]
\thicklines
\centering{
\begin{picture}(25,90)(-18,30)

\put(-123,110){\small{$\boldsymbol{W}$}}
\put(36,110){\small{$\boldsymbol{X}$}}
\put(147,110){\small{$\boldsymbol{Y}$}}

\put(36,82){\small{$\boldsymbol{P}$}}


\put(-160,40){\circle*{3}}
\put(-120,40){\circle*{3}}

\put(-160,40){\line(-1,4){15}}
\put(-160,40){\line(0,1){60}}
\put(-160,40){\line(1,4){15}}

\put(-168,70){\circle*{3}}
\put(-160,70){\circle*{3}}
\put(-152,70){\circle*{3}}

\put(-128,70){\circle*{3}}
\put(-120,70){\circle*{3}}
\put(-112,70){\circle*{3}}

\put(-88,70){\circle*{3}}
\put(-80,70){\circle*{3}}
\put(-72,70){\circle*{3}}

\put(-80,40){\circle*{3}}

\put(-120,40){\line(-1,4){15}}
\put(-120,40){\line(0,1){60}}
\put(-120,40){\line(1,4){15}}

\put(-80,40){\line(-1,4){15}}
\put(-80,40){\line(0,1){60}}
\put(-80,40){\line(1,4){15}}

\put(-181,67){\small{$a_x$}}
\put(-141,67){\small{$a_y$}}
\put(-101,67){\small{$a_z$}}

\put(1,67){\small{$a$}}

\put(-163,30){\small{$x$}}
\put(-123,30){\small{$y$}}
\put(-83,30){\small{$z$}}


\put(-34,84){\small{$g_1$}}

\put(-40,77){\vector(1,0){20}}
\put(-40,71){\vector(1,0){20}}

\put(-34,61){\small{$g_2$}}


\put(10,40){\circle*{3}}
\put(70,40){\circle*{3}}
\put(10,40){\line(0,1){30}}
\put(10,70){\circle*{3}}

\put(20,70){\small{$\dots$}}
\put(49,70){\small{$\dots$}}

\put(10,75){\line(0,1){3}}
\put(10,81){\line(0,1){3}}
\put(10,87){\line(0,1){3}}
\put(10,93){\line(0,1){3}}

\put(70,75){\line(0,1){3}}
\put(70,81){\line(0,1){3}}
\put(70,87){\line(0,1){3}}
\put(70,93){\line(0,1){3}}

\put(10,40){\line(1,1){30}}
\put(10,40){\line(2,1){60}}

\put(70,40){\line(-1,1){30}}
\put(70,40){\line(-2,1){60}}

\put(40,70){\circle*{3}}
\put(70,40){\line(0,1){30}}
\put(70,70){\circle*{3}}
\put(7,30){\small{$x$}}
\put(67,30){\small{$y$}}


\put(101,77){\small{$h$}}

\put(94,71){\vector(1,0){20}}


\put(142,77){\small{$h(x)$}}
\put(137,59){\small{$=\!h(y)$}}
\put(150,70){\circle*{3}}
\put(150,70){\circle{50}}

\end{picture}}\nopagebreak
\caption{} \label{fig:1}
\end{figure}
}

This delivers the desired contradiction, because $\sbW\in\mathsf{K}_*$.  Indeed, $\mathsf{K}$ is a variety containing $\sbX^*$ and, by
Lemma~\ref{easy}(\ref{upspace lem}),(\ref{disjoint union}),
$(\bup{\sbx})^*\in\mathbb{H}(\sbX^*)$ and
$\sbW^*\in\mathbb{IP}((\bup{\sbx})^*)$, so $\sbW^*\in\mathsf{K}$.

\smallskip

{\bf Case~II.}  \,$x$ has depth $n$ and $y$ has depth less than $n$.

\smallskip

In this case, $n\geq 2$ and $x<y$, by Lemma~\ref{claim b}.  In fact, $y$ covers $x$, because $h|_P$ is injective.  Moreover, $y$ is the only cover of $x$.
(For, if $u$ is another cover, then $h(y)=h(x)\leqslant h(u)$, whence $h(u)=h(v)$ for some $v\in\up{y}$.  But then, $u$ and $v$ are distinct
and have depth less than $n$, contradicting the injectivity of $h|_P$.)  Let $\bup{\sbz}$ (with least element $z$) be a disjoint copy of $\bup{\sbx}$
($=\sbX$),
and $\sbW$ the disjoint union of $\bup{\sbx}$ and $\bup{\sbz}$.  As in Case~I, \,$\sbW\in\mathsf{K}_*$ and we can construct distinct
Esakia morphisms $g_1,g_2\colon\sbW\mrig\sbX$ with $h\circ g_1=h\circ g_2$ (where $g_1(z)=x$ and $g_2(z)=y$).
\end{proof}

In a variety $\mathsf{K}$, when we verify that a $\mathsf{K}$--morphism $h\colon\sbA\mrig\sbB$ is a $\mathsf{K}$--epimorphism,
it is enough to show that, for any homomorphisms $f,g$ from $\sbB$ to a {\em subdirectly irreducible\/} member of $\mathsf{K}$,
if $f\circ h= g\circ h$, then $f=g$.  This is simply a consequence of the subdirect decomposition theorem.

Note also that a variety $\mathsf{K}$ has the ES property iff all epimorphic {\em inclusion\/} maps in $\mathsf{K}$ are surjective,
because a $\mathsf{K}$--morphism $h\colon\sbA\mrig\sbB$ is a $\mathsf{K}$--epimorphism iff
the inclusion $h[\sbA]\mrig\sbB$ is.  Thus, $\mathsf{K}$ has the ES property iff
no $\sbB\in\mathsf{K}$ has a proper subalgebra $\sbC$ that is $\mathsf{K}$--{\em epic\/} in the sense that $\mathsf{K}$--morphisms with domain $\sbB$
are determined by their restrictions to $\sbC$.

\begin{thm}\label{main brouwerian}
Let\/ $\mathsf{K}$ be a variety of Brouwerian algebras, where\/ $\mathsf{K}$ has finite depth.  Then epimorphisms in\/ $\mathsf{K}$ are
surjective.
\end{thm}
\begin{proof}
Let $\mathsf{K}_\bot$ be the subvariety of $\mathsf{HA}$ generated by
$\{\sbA_\bot:\sbA\in\mathsf{K}\}$.
By Theorems~\ref{heyting depth formulas}--\ref{brouwerian depth formulas} and \ref{new main},
$\mathsf{K}_\bot$ has finite depth, and hence the ES property.
Suppose $\sbA$ is a $\mathsf{K}$--epic subalgebra of some $\sbB\in\mathsf{K}$.  We must show that $A=B$.

Identifying $\sbA_\bot$ appropriately with a subalgebra of $\sbB_\bot$, we claim that it is a
$\mathsf{K}_\bot$--epic subalgebra.  To see this,
let $f,g\colon\sbB_\bot\mrig\sbD$ be $\mathsf{K}_\bot$--morphisms that agree on $\sbA_\bot$.
Recall that $\sbD$ may be assumed subdirectly irreducible, in which case
$\sbD=\sbE_\bot$ for some $\sbE\in\mathsf{K}$, by Lemma~\ref{another bot lem}.
Then $f|_B$ and $g|_B$ are $\mathsf{K}$--morphisms $\sbB\mrig\sbE$ (Remark~\ref{bot congruences}),
which agree on $\sbA$, so $f|_B=g|_B$, as $\sbA$ is $\mathsf{K}$--epic in $\sbB$. This forces $f=g$, so
$\sbA_\bot$ is indeed $\mathsf{K}_\bot$--epic in $\sbB_\bot$.  Thus,
$A_\bot=B_\bot$, by the ES property of $\mathsf{K}_\bot$, and so
$A=B$.
\end{proof}

\begin{cor}\label{main}
Epimorphisms are surjective in every finitely generated variety of Heyting or Brouwerian algebras.
\end{cor}
\begin{proof}
Use Theorems~\ref{new main} and \ref{main brouwerian}
and Lemma~\ref{finite}.
\end{proof}

An argument of Kuznetsov \cite{Kuz75} shows that $2^{\aleph_0}$ varieties of Heyting algebras
(and as many of Brouwerian algebras) have depth $3$.  So, among the subvarieties of $\mathsf{HA}$
or of $\mathsf{BRA}$, a continuum
have the ES property, by Theorems~\ref{new main} and \ref{main brouwerian}.
Only denumerably many of these are finitely generated, and only finitely many have the strong ES property \cite{Mak00,Mak03}.

A logic algebraized by a finitely generated or locally finite variety is said to be {\em tabular\/}
or {\em locally tabular}, respectively.

\begin{cor}\label{int tabular inf beth}
If a super-intuitionistic or positive super-intuitionistic logic is tabular---or more generally, if its theorems include\/ $h_n$
for some\/ $\textup{$n\in\omega$}$---then it
has the infinite Beth property.
\end{cor}

A subdirect product of totally ordered Heyting or Brouwerian algebras is called a {\em G\"{o}del algebra\/} or a {\em relative Stone algebra},
respectively.  These form varieties $\mathsf{GA}$ and $\mathsf{RSA}$, whose respective subvarieties algebraize
the {\em G\"{o}del logics\/} and the {\em positive G\"{o}del logics}, cf.\ \cite{Dum59,Haj98}.
It is well known (and implicit in \cite{DW73}) that every subquasivariety of $\mathsf{GA}$
or of $\mathsf{RSA}$ is a variety.  G\"{o}del algebras are examples of the {\em BL-algebras\/} of \cite{Haj98}, which algebraize
{\em $\text{Hajek's}$ basic logic}.  The subvarieties of $\mathsf{GA}$ are the only varieties of BL-algebras having the weak ES property \cite{Mon06}.

\begin{cor}\label{godel}
Epimorphisms are surjective in every variety of G\"{o}del algebras or of relative Stone algebras.  In other words, all
G\"{o}del logics and positive G\"{o}del logics have the infinite Beth property.
\end{cor}
\begin{proof}
This follows from Corollary~\ref{main}, because the
only subvarieties of $\mathsf{GA}$ or $\mathsf{RSA}$ that are not finitely generated are $\mathsf{GA}$ and $\mathsf{RSA}$
themselves \cite{DM71}, and they
have the strong ES property.
\end{proof}
\noindent
Indeed, $\mathsf{GA}$ and $\mathsf{RSA}$ both have denumerably many subvarieties, and the ones with the strong ES property
are $\mathsf{GA}$,
$\mathbb{V}(\sbG_3)$, the variety of Boolean algebras, the trivial variety of Heyting algebras, and the Brouwerian subreduct
classes of these \cite{Mak00,Mak03}.  Here, $\sbG_3$ denotes the three-element G\"{o}del algebra.

$\mathsf{GA}$ and $\mathsf{RSA}$ have infinite depth, but they are the largest subvarieties of $\mathsf{HA}$ and $\mathsf{BRA}$
(respectively) having width $1$.  As we shall see in the next section, width $2$ is not a sufficient condition for the ES property.

\section{A Counter-Example}\label{a counter-example section}

In this section, we establish Blok and Hoogland's conjecture that the ES and weak ES properties are distinct, by exhibiting a
variety of Brouwerian algebras (and a variety of Heyting algebras) in which not all epimorphisms are surjective.  Recall that this is sufficient,
because all varieties of Brouwerian or Heyting algebras have the weak ES property, by Theorem~\ref{heyting es}(\ref{heyting es1}).
As it happens, the counter-examples can be chosen locally finite, so Corollary~\ref{main} cannot be generalized to all locally finite
varieties.

Let $\A$ be the denumerable subdirectly irreducible Brouwerian algebra whose lattice reduct is depicted in Figure~\ref{fig:2}.

{\small
\begin{figure}[htbp]
\thicklines
\centering{
\begin{picture}(250,163)(0,12)

\put(125,138){\line(0,1){15}}
\put(125,153){\circle*{3}}

\put(125,108){\circle*{3}}
\put(140,123){\line(-1,-1){15}}
\put(140,123){\circle*{3}}
\put(110,123){\line(1,-1){15}}
\put(110,123){\circle*{3}}
\put(110,123){\line(1,1){15}}
\put(125,138){\circle*{3}}
\put(140,123){\line(-1,1){15}}

\put(125,78){\circle*{3}}
\put(140,93){\line(-1,-1){15}}
\put(140,93){\circle*{3}}
\put(110,93){\line(1,-1){15}}
\put(110,93){\circle*{3}}
\put(110,93){\line(1,1){15}}
\put(125,108){\circle*{3}}
\put(140,93){\line(-1,1){15}}

\put(125,48){\circle*{3}}
\put(140,63){\line(-1,-1){15}}
\put(140,63){\circle*{3}}
\put(110,63){\line(1,-1){15}}
\put(110,63){\circle*{3}}
\put(110,63){\line(1,1){15}}
\put(125,78){\circle*{3}}
\put(140,63){\line(-1,1){15}}

\put(140,63){\line(-1,-1){22}}
\put(110,63){\line(1,-1){22}}

\put(125,30){\line(0,-1){2}}
\put(125,26){\line(0,-1){2}}
\put(125,22){\line(0,-1){2}}

\put(122,160){$\boldmath{\top}$}

\put(97,122){\small $b_0$}
\put(146.5,122){\small $c_0$}

\put(97,92){\small $b_1$}
\put(146.5,92){\small $c_1$}

\put(97,62){\small $b_2$}
\put(146.5,62){\small $c_2$}

\end{picture}}\nopagebreak
\caption{} \label{fig:2}
\end{figure}
}

In a partially ordered set, two elements will be called {\em siblings\/} if they are incomparable.
Note that $\sbA$ obeys the {\em sibling rule\/}: each of its elements has at most one sibling.

\begin{Theorem}\label{Thm : ES fails}
The ES property fails in\/ $\VVV(\A)$ and in\/ $\VVV(\A_{\bot})$.
\end{Theorem}
\begin{proof}
The sibling rule is expressed
by the following positive universal sentence in the language of $\mathsf{BRA}$
(where $x \leqslant y$ abbreviates $x \land y \thickapprox x$):
\begin{equation}\label{antichain}
\forall x\,\forall y\,\forall z \left( x \leqslant y \, \textsf{ or } \, y \leqslant x \, \textsf{ or } \, x \leqslant z \, \textsf{ or } \, z \leqslant x \, \textsf{ or } \, y \thickapprox z\right)\!.
\end{equation}
Positive universal sentences persist under
$\mathbb{H}$, $\mathbb{S}$ and $\mathbb{P}_\mathbb{U}$, so by J\'{o}nsson's Theorem,
every
member of $\mathbb{V}(\sbA)_\textup{FSI}$ obeys the sibling rule.
(It is then easily verified that $\mathbb{V}(\sbA)$
has width $2$, and likewise $\mathbb{V}(\sbA_\bot)$.)

Let $\B$ be the subalgebra of $\A$ with universe
$B=\{ b_{i} : i \in \omega \} \cup \{ \top \}$.
As $B\neq A$, the ES property for $\VVV(\sbA)$ will be refuted if we can show that $\sbB$ is $\VVV(\sbA)$--epic in $\sbA$.

Suppose not.
Then there exist $\C \in \VVV(\A)$ and distinct homomorphisms $f, g \colon \A \mrig \C$ that agree on $B$.  Moreover,
$\sbC$ may be assumed subdirectly irreducible, so it obeys the sibling rule.

Consider the filters $F\seteq f^{-1}[\{\top\}]$ and $G=g^{-1}[\{\top\}]$ of $\sbA$.

As $A=\up{B}$ and $B\cap F=B\cap G$, and since $F$ and $G$ are up-sets of $\langle A;\leqslant\rangle$,
neither $F$ nor $G$ can be $A$.  (Otherwise, $f$ and $g$ would both be the constant function with range $\{\top\}$.)
Because $\sbC$ is subdirectly irreducible, its subalgebras are finitely subdirectly irreducible, so
$F$ and $G$ are prime.
\begin{Claim}\label{claim 1}
There exists $i \in \omega$ such that $b_{i}$ is the least element of $B \cap F$.
\end{Claim}

Indeed, since $F\neq A$, the Hasse diagram shows that $B\cap F$ has a least element, $a$ say.  It remains to verify that $a\neq\top$.  Suppose $a=\top$.
Then
\begin{equation}\label{new eq}
F, G \subseteq \{ c_{0}, b_{0} \lor c_{0}, \top \}.
\end{equation}
Observe that
\begin{equation}\label{Eq : Incomparable elements}
b_{i} \to c_{i} = c_{i} \text{ \,and\, }c_{i} \to b_{i} = b_{i}\text{ \,for every } i \in \omega.
\end{equation}
Consider any $i \geq 1$. By (\ref{Eq : Incomparable elements}), $b_{i} \to c_{i},\, c_{i} \to b_{i} \notin F \cup G$,
so $f(b_i)\rig f(c_i)$ and $f(c_i)\rig f(b_i)$ are not $\top$, i.e.,
$f(b_{i})$ and $f(c_{i})$ are siblings, and likewise $g(b_{i})$ and $g(c_{i})$.
But $f(b_{i}) = g(b_{i})$, so the sibling rule forces
$f(c_{i}) = g(c_{i})$. Applying joins, we see that
$f(a) = g(a)$ for every $a \leqslant b_{0}$, so $f$ and $g$ must disagree on $X \seteq \{ c_{0}, b_{0} \lor c_{0} \}$.
This rules out the possibility that $f[X]=\{\top\}=g[X]$, so
$F$ and $G$ can't both be
$\up{c_{0}}$.
If $c_0\in F$, then $F=\up{c_0}$, by (\ref{new eq}), and
$f(b_{0}) = f(b_{0}\land c_{0})$, but $f$ and $g$ agree on $\{ b_{0}, b_{0} \land c_{0} \}$, so $c_{0} = b_{0} \to (b_{0} \land c_{0}) \in G$,
whence $G=\up{c_0}=F$, contradicting the previous sentence.  Therefore, $c_0\notin F$ and, by symmetry, $c_0\notin G$.
Then $f(b_{0})$ and $f(c_{0})$ are siblings, by (\ref{Eq : Incomparable elements}), and likewise the pair $g(b_{0}),g(c_{0})$.
As $f(b_0)=g(b_0)$, the sibling rule yields $f(c_{0}) = g(c_{0})$, whence
$f(b_{0} \lor c_{0}) = g(b_{0} \lor c_{0})$.  But then, $f = g$, a contradiction.   This vindicates Claim~\ref{claim 1}.

As $F$ and $G$ are prime and $B\cap F=B\cap G$, it follows from Claim~\ref{claim 1} and the Hasse diagram that
$F \subseteq G$ or $G \subseteq F$.  By symmetry, we may assume that $G\subseteq F$.
There are then two cases: $G\subsetneq F$ or $F=G$.  In each of these cases, we shall obtain a contradiction, as desired.

\smallskip

{\bf Case~I.}  \ $G\subsetneq F$.

\smallskip

By the primeness of $F$ and $G$,
\[
F = \up{c_{i+1}} \text{ \,and\, }G = \up{b_{i}} \textup{ \ (for the $i$ in Claim~\ref{claim 1})}.
\]
Let $j \geq i+2$.  Then $b_{j} \to c_{j},\, c_{j} \to b_{j} \notin F \cup G$.  So, applying the sibling rule (as above), we infer from $f(b_j)=g(b_j)$
that $f(c_{j}) = g(c_{j})$.  Then, applying joins, we obtain $f(a) = g(a)$ for every $a \leqslant b_{i+1}$.
Now observe that
\[
\{ b_{i+1} \to ( b_{i+1} \land c_{i+1} ),\,\, (b_{i+1} \land c_{i+1} )\to b_{i+1} \}
\]
is a subset of $F$, but not of $G$.  In other words, $f(b_{i+1}) = f(b_{i+1} \land c_{i+1})$ but $g(b_{i+1}) \ne g(b_{i+1} \land c_{i+1})$, contradicting the fact that $f$ and $g$ agree on $\down{b_{i+1}}$.

\smallskip

{\bf Case~II.}  \ $F=G$.

\smallskip

For the $i$ in Claim~\ref{claim 1}, we have $F = G = \up{b_{i}}$ or $F = G = \up{c_{i+1}}$.  We deal only with the former case, because the latter is analogous.  Observe that
$b_{j} \to c_{j},\, c_{j} \to b_{j} \notin F $ for every $j \geq i+1$.  So, because $f$ and $g$ agree on $B$, the sibling rule yields $f(a) = g(a)$ for
every $a \leqslant b_{i} \land c_{i}$.  On the other hand, $f(a) = \top = g(a)$ for every $a \geqslant b_{i}$. Thus, $f$ and $g$ coincide on $A\smallsetminus\{c_{i}\}$.
But $c_{i} \to (b_{i} \land c_{i}),\, (b_{i} \land c_{i}) \to c_{i} \in F$, so $f(c_{i}) = f(b_{i} \land c_{i}) = g(b_{i} \land c_{i}) = g(c_{i})$, whence $f= g$, a contradiction.

We have shown that $\sbB$ is $\VVV(\sbA)$--epic in $\sbA$, whence
$\VVV(\A)$ lacks the ES property.
With only notational changes,
the same argument shows that the subalgebra $\B_{\bot}$ of $\sbA_\bot$ is $\VVV(\A_{\bot})$--epic,
so the ES property fails for $\VVV(\A_{\bot})$ too.  (Here,
$\bot\notin F\cup G$,
because $\sbC$ is nontrivial.)
\end{proof}

In contrast, a finite subalgebra $\sbD$ of $\sbA_\bot$ generates a variety with the ES property (Corollary~\ref{main}), so
$\sbB_\bot\cap\sbD$ cannot be a proper $\mathbb{V}(\sbD)$--epic sub\-algebra of $\sbD$.  For example, the subalgebra of $\sbA_\bot$
generated by $\{b_0,c_0\}$ has two endomorphisms sending $b_0$ to $\top$ but disagreeing at $c_0$.  (The pre-image of $\{\top\}$ is
$\up{b_0}$ in one case, and $\up{(b_0\wedge c_0)}$ in the other.)

In \cite{Esa85,Esa85b}, it is shown that the variety of Heyting algebras satisfying the {\em weak Peirce law\/}
\begin{equation}\label{weak peirce}
(y \to x) \lor (((x \to y) \to x) \to x) \thickapprox \top
\end{equation}
is locally finite.
It follows that the Brouwerian algebras satisfying (\ref{weak peirce}) also form a locally finite variety, because if $\sbE\in\mathsf{BRA}$ is
$n$--generated and satisfies (\ref{weak peirce}),
then the same is true of
$\sbE_\bot$.
It is easily checked that $\A$ satisfies (\ref{weak peirce}), so $\VVV(\A_{\bot})$ and $\VVV(\sbA)$ are both locally finite, and we have proved:

\begin{Corollary}\label{blok hoogland answer}
The ES and weak ES properties are distinct, even for locally finite varieties of Brouwerian or Heyting algebras having width\/ $2$\textup{.}
\end{Corollary}

A subvariety $\mathsf{M}$ of $\mathsf{BRA}$ or $\mathsf{HA}$ satisfies
an identity of the form
\[
\al_1\vee\ldots\vee\al_n\thickapprox\top
\]
iff the formula
$\al_1\thickapprox\top\textup{ \,$\mathsf{or}$\, } \dots \textup{ \,$\mathsf{or}$\, }
\al_n\thickapprox\top$
is valid in $\mathsf{M}_\textup{FSI}$ (because $\mathsf{M}_\textup{FSI}$ comprises the algebras in $\mathsf{M}$ where $\top$ is join-irreducible).
The variety generated by the Brouwerian algebras satisfying the sibling rule (\ref{antichain}) is therefore axiomatized,
relative to $\mathsf{BRA}$, by
\begin{equation}\label{Eq : Axioms}
(x \to y) \lor (y \to x) \lor (x \to z) \lor (z \to x) \lor (y \leftrightarrow z) \thickapprox \top,
\end{equation}
where $y\leftrightarrow z$ abbreviates $(y\rig z) \wedge (z\rig y)$.

Clearly, this variety contains $\mathbb{V}(\sbA)$, but in fact
they can be shown equal,
i.e., a Brouwerian algebra satisfies (\ref{Eq : Axioms}) iff it belongs to $\mathbb{V}(\sbA)$.
On the other hand, $\VVV(\sbA_\bot)$ is axiomatized, relative to $\mathsf{HA}$, by (\ref{Eq : Axioms}) and
the {\em weak excluded middle\/} identity
$\lnot x \lor \lnot \lnot x \thickapprox \top$.
Proofs of these claims will
appear in a subsequent paper,
where we shall also construct $2^{\aleph_0}$ sub\-varieties of $\mathsf{HA}$ without the ES property.

Corollary~\ref{blok hoogland answer} and the correspondences in Section~\ref{introduction} show that the finite and infinite Beth properties
are distinct, even in the context of
locally tabular super-intuitionistic logics.

Specifically, let $\mathbf{L}$ be the extension
of $\mathbf{IPL}$ by the axioms
\[
\lnot x \lor \lnot \lnot x \quad \textup{and} \quad
(x \to y) \lor (y \to x) \lor (x \to z) \lor (z \to x) \lor (y \leftrightarrow z),
\]
so $\mathbf{L}$ (more exactly, $\,\vdash_\mathbf{L}$) is algebraized by $\VVV(\sbA_\bot)$.
A failure of the infinite Beth property
can be extracted from the $\mathbb{V}(\sbA_\bot)$--epimorphic inclusion $i\colon\sbB_\bot\mrig\sbA_\bot$ by a general method
in \cite{BH06},
which can be made more concrete for $\mathbf{L}$.  We create disjoint sets of distinct variables
\[
X = \{ x_{a} : a \in B_{\bot} \} \text{ \,and\, }Z = \{ z_{a} : a \in A_{\bot} \smallsetminus B_{\bot} \},
\]
and with each connective $\ast \in \{ \land, \lor, \to \}$, we associate a set of formulas
\[
\Sigma_{\ast} \seteq \{ (u \ast v) \leftrightarrow w : \,u, v, w \in X \cup Z \text{ \,and\, } f(u) \ast f(v) = f(w) \},
\]
where $f\colon X\cup Z\mrig A_\bot$ sends each variable to its subscript.  Then
\[
\Gamma \seteq \,\Sigma_{\land} \,\cup\, \Sigma_{\lor} \,\cup\, \Sigma_{\to} \,\cup\, \{ x_\top \leftrightarrow \top,\,x_\bot \leftrightarrow \bot \}
\]
captures the infinite `diagram' of $\sbA_\bot$.
For an $\mathbf{L}$--formula $\psi$ over $X\cup Z$,
\begin{equation}\label{theory}
\textup{$\Gamma\vdash_\mathbf{L}\psi$ iff $\psi$ takes the value $\top$ when $f$ interprets $X\cup Z$
in $\sbA_\bot$.}
\end{equation}
The forward implication follows by induction on the length of a proof of $\psi$ from $\Gamma$ in $\mathbf{L}$, and its converse from the fact
that
$\Gamma\vdash_\mathbf{L}\psi\leftrightarrow v$ for some $v\in X\cup Z$ (which in turn follows by induction on the complexity of $\psi$).

As $i$ is a $\mathbb{V}(\sbA_\bot)$--epimorphism, (\ref{theory}) and
the proof of \cite[Thm.\,3.12]{BH06} show that
$\Gamma$ defines $Z$ implicitly in terms of $X$ in $\mathbf{L}$, i.e., $\Gamma\cup \sigma
[\Gamma]\vdash_\mathbf{L} z\leftrightarrow \sigma(z)
$
for all $z\in Z$ and all substitutions $\sigma$
(defined on $X\cup Z$) that fix every element of $X$.
Using $f$ to evaluate $X\cup Z$ in $\sbA_\bot$, we see
that $Z$ is not defined
explicitly in terms of $X$ by $\Gamma$ in $\mathbf{L}$, because a formula in $\Gamma$, a formula
over $X$ and an element of $Z$ take values in $\{\top\}$, $B_\bot$ and $A_\bot\smallsetminus B_\bot$, respectively.

A more economical refutation of the infinite Beth property in $\mathbf{L}$ can be inferred from this example.  Writing $x_i$ and $z_i$
for $x_{b_i}$ and $z_{c_i}$, respectively, let $X'=\{x_i:i\in\omega\}$, $Z'=\{z_i:i\in\omega\}$ and
\begin{align*}
& \Gamma'\,=\,\,\mbox{$\bigcup$}_{\,i\,\in\,\omega}\,\{(x_i\wedge z_i)\leftrightarrow (x_{i+1}\vee z_{i+1}),\,\neg\neg x_i,\,\neg\neg z_i\}\\
& \quad\,\cup\;\,\mbox{$\bigcup$}_{\,j\,>\,i\,}\,\{(x_i\rig x_j)\rig x_j,\,(z_i\rig z_j)\rig z_j\}\\
& \quad\,\cup\;\,\mbox{$\bigcup$}_{\,j\,\geq\, i\,}\,\{(x_i\rig z_j)\rig z_j,\,(z_i\rig x_j)\rig x_j\}.
\end{align*}
Let $\rho$ be the substitution over $X\cup Z$ that fixes each element of $X'\cup Z'$ while sending $x_\top$ to $\top$, $x_\bot$ to $\bot$,
and each $z_{b_i\vee c_i}$ to $x_i\vee z_i$.

On algebraic grounds, $\Gamma'\vdash_{\mathbf{IPL}}\rho[\Gamma]$,
so it follows readily from the previous example that $\Gamma'$ defines $Z'$ implicitly in terms of $X'$ in $\mathbf{L}$, and not explicitly.
In other words, no formulas over $X'$, constrained by the relations expressed in $\Gamma'$, define values for the variables $z_i$ in all
algebraic models of $\mathbf{L}$, yet any interpretation of $X'$ in such a model uniquely determines (or precludes) values for each $z_i$,
subject to the same relations.

\section{Non-Integral Varieties}\label{non-integral varieties section}

Recall that a category equivalence $F$ between two varieties
induces an isomorphism $\overline{F}\colon\mathsf{M}\mapsto\mathbb{I}\{F(\sbA)\colon\sbA\in\mathsf{M}\}$ between their subvariety lattices.
Moreover,
$F$ restricts to a category equivalence from
$\mathsf{M}$
to
$\overline{F}(\mathsf{M})$
for each subvariety $\mathsf{M}$.
In this situation, $\mathsf{M}$ is finitely generated iff $\overline{F}(\mathsf{M})$ is.
The ES property and its weak and strong analogues are preserved by any category equivalence
between varieties.
Another invariant
is the demand that all subquasivarieties be varieties.  (All of these
claims are justified in \cite[pp.\,222 \& 238]{McK96} and/or \cite[Sec.\,5 \& 7]{GR12}, for instance.)

\begin{exmp}
{\em
It follows from Corollary~\ref{godel} that all axiomatic extensions of the logics called $\mathbf{IUML}$ and $\mathbf{IUML}^-$ in \cite{MM10,MM07}
have the infinite Beth property, as the algebraic counterparts of these two systems are categorically equivalent to
$\mathsf{GA}$ and $\mathsf{RSA}$ \cite{GR12}, and the ES property for all subvarieties is transferred by the equivalence.
The latter case instantiates a wider ES result (Theorem~\ref{sugihara es}) in the next section,
so we
postpone further discussion of the algebras---except to say that they are CRLs which need not be integral.
}\qed
\end{exmp}

The paper \cite{GR15} establishes further category equivalences facilitating the transfer of information from integral
to non-integral settings (the former being better understood at present).  For the fullest exploitation of
this strategy in the case of ES properties, the next theorem is helpful.

\begin{thm}\label{expansions}
Let\/ $\mathsf{K}$ be a variety consisting of expansions of Brouwerian algebras, where the additional operations on each\/ $\sbA\in\mathsf{K}$ are
compatible with the congruences of the Brouwerian reduct\/ $\sbA^-$ of\/ $\sbA$\textup{.}  Then
\begin{enumerate}
\item\label{expansions 1}
$\mathsf{K}$ has the weak ES property, and

\item\label{expansions 2}
if\/ $\mathsf{K}$ is finitely generated, it
has the ES property.
\end{enumerate}
\end{thm}
\begin{proof}
(\ref{expansions 1}) follows from \cite[Thm.\,12.3]{GR15}, as the weak ES property is the algebraic counterpart of the finite Beth property.

(\ref{expansions 2})
Let $\sbC$ be a finite algebra such that $\mathsf{K}=\mathbb{V}(\sbC)$.  Let
$i\colon\sbA\mrig\sbB$ be an inclusion map that is a $\mathsf{K}$--epimorphism.  It suffices to show that $i$ is surjective
(i.e., that $A=B$).  Suppose not.

Viewing $i$ as a $\mathsf{BRA}$--morphism $\sbA^-\mrig\sbB^-$, we infer from
Lemma~\ref{easy}(\ref{easy 1}) that the dual map $i_*\colon\sbB^-_*\mrig\sbA^-_*$ is not injective,
i.e., there are distinct $x,y\in
\textup{Pr}\,\sbB^-$ with $i_*(x)=i_*(y)$.

In $\sbB^-_*$, the closed up-set $W\seteq\up{\{x,y\}}$
is the universe of an E-subspace $\sbW$.
Let $j\colon\sbW\mrig\sbB^-_*$ be the inclusion map.
Then $\sbW^*\in\mathbb{H}(\sbB^-)$, by Lemma~\ref{easy}(\ref{upspace lem}), so
we may assume that $\sbW^*=\sbB^-/\theta$ for some $\theta\in{\boldsymbol{\mathit{Con}}}\,\sbB^-$
and that $j^*$ is the canonical map $\lambda_\theta\colon b\mapsto b/\theta$.  By the assumption in the present theorem's statement, $\theta\in{\boldsymbol{\mathit{Con}}}\,\sbB$,
so $\lambda_\theta\colon\sbB\mrig\sbB/\theta$ is a $\mathsf{K}$--morphism, hence a $\mathsf{K}$--epimorphism (being surjective).
Then $\lambda_\theta\circ i$
is also a $\mathsf{K}$--epimorphism, because $i$ is.

Now $\sbB^-\in\mathbb{V}(\sbC^-)$,
which is a finitely generated subvariety of $\mathsf{BRA}$, so $W$ is finite, by Lemma~\ref{finite}.  Then $\sbB/\theta$ is finite,
because $\sbB^-/\theta$ is the dual of $\sbW$.  Therefore, by (\ref{expansions 1})
and the remark before Definition~\ref{weak es def},
the $\mathsf{K}$--epimorphism $\lambda_\theta\circ i\colon\sbA\mrig\sbB/\theta$ is
surjective.
Its dual $i_*\circ j$ must then be injective, by Lemma~\ref{easy}(\ref{easy 1}).  But this is a contradiction, because $x,y$ belong
to the domain $W$ of $i_*\circ j$, and $i_*j(x)=i_*(x)=i_*(y)=i_*j(y)$, while $x\neq y$.
\end{proof}

For example, a \textit{nuclear relative Stone algebra\/} $\sbA$ is the expansion of a relative Stone algebra by a \textit{nucleus}, i.e., by a unary operation $\lozenge$ such
that
\[
\textup{$a\leqslant\lozenge a=\lozenge\lozenge a$ \,and\,
$\lozenge a \wedge \lozenge b=\lozenge (a\wedge b)$}
\]
for all $a,b\in A$.
Any such algebra has the same congruences as its Brouwerian reduct
\cite[Thm.\,7.1]{GR15}, so Theorem~\ref{expansions} applies.  Although algebras of this kind have been studied independently,
they serve here as stepping stones to non-integral residuated structures (without nuclei).  The latter model logics that lack
the {\em weakening axiom\/} $x\rig(y\rig x)$.

Specifically, by \cite[Cor.\,3.5]{GR15}, a subdirect product of totally ordered idempotent (possibly non-integral) CRLs is generated by the lower bounds of its neutral element $\sbe$ iff it satisfies
$((x\vee \sbe)\rig \sbe)\rig \sbe\,\thickapprox\,x\vee \sbe$.
In this case, it is called a {\em generalized Sugihara monoid}.  These algebras form a variety, which is shown in \cite[Thm.\,8.7]{GR15} to be
categorically equivalent to the variety of nuclear relative Stone algebras.
The category equivalence adapts to the corresponding bounded cases \cite[p.\,3208]{GR15}.
So, by Theorem~\ref{expansions}(\ref{expansions 1}), every variety of generalized Sugihara monoids or bounded ones has the
weak ES property.  That was already observed in \cite[Thm.\,13.1]{GR15}, but Theorem~\ref{expansions}(\ref{expansions 2}) and the opening paragraph of this
section also yield:

\begin{thm}\label{gsm}
Epimorphisms are surjective in every finitely generated variety consisting of
generalized Sugihara monoids or bounded ones.
\end{thm}

When {\em negation\/} connectives of substructural logics are modeled in non-integral CRLs, they can normally be identified with functions
$\textup{$a\mapsto a\rig \sbf$}$, where $\sbf$ is a fixed element of the algebra, but not the least element.  That deprives $\bot$ of
its \textit{raison d'\^{e}tre} from
intuitionistic logic, so it is often discarded from the signature.  It is partly for this reason that we have been attentive to Brouwerian
(not only Heyting) algebras thus far.
In particular, bounds are traditionally neglected
in the residuated structures
considered below.

\section{Sugihara Monoids}\label{sugihara monoids section}

A unary operation $\sim$ on a CRL $\boldsymbol{A}$
will be called an
\textit{involution} if
\[
\textup{$\simop \simop a = a$ \,and\,
$a \to \simop b = b \to \simop a$}
\]
for all $a, b \in A$ (in which case, $\simop a=a\rig\simop \sbe$ for all $a$).
The
expansion of $\boldsymbol{A}$ by $\sim$ is then called an \textit{involutive CRL}.

The variety $\mathsf{SM}$ of \textit{Sugihara monoids} comprises the idempotent distributive involutive CRLs.
It algebraizes the relevance logic $\mathbf{RM}^\mathbf{t}$ of \cite{AB75}.  Sugihara monoids are not integral,
unless they are Boolean algebras.

We denote by $\mathsf{X}$ the class of nuclear relative Stone algebras with a distinguished element $\boldsymbol{f}$ satisfying the following quasi-equations:
\begin{eqnarray*}
& x \lor (x \to \boldsymbol{f}) \,\thickapprox\, \top \,\thickapprox\, \lozenge ( \lozenge x \to x )\\
& \lozenge x \thickapprox \top \,\Longleftrightarrow\, \boldsymbol{f} \leqslant x.
\end{eqnarray*}
It can be proved that $\mathsf{X}$ is a variety (see \cite[p.\,3207]{GR15}). We shall need:

\begin{Theorem}\label{Theorem : equivalence}
\textup{(\cite[Thm.\,10.5]{GR15})}
\,$\mathsf{SM}$ and\/ $\mathsf{X}$ are categorically equivalent.
\end{Theorem}
Not all proper subvarieties of $\mathsf{SM}$ are finitely generated, so Theorem~\ref{expansions}(\ref{expansions 2}) has limited utility for Sugihara monoids.
Nevertheless, we aim to establish the ES property for {\em all\/} subvarieties of $\mathsf{SM}$.  By
Theorem~\ref{Theorem : equivalence},
it suffices to prove the corresponding result for the simpler variety $\mathsf{X}$, where we can exploit Esakia duality in the Brouwerian reducts of the algebras.
The next lemma is just a specialization of Theorem~\ref{expansions}(\ref{expansions 1}).

\begin{Lemma}\label{Lemma : X has weak ES}
Every subvariety of\/ $\mathsf{X}$ has the weak ES property.
\end{Lemma}

In \cite[Thms.\,\,4.4 \& 4.9]{MM12}, the subvariety lattice of $\mathsf{SM}$ is described, and the amalgamable subvarieties are identified.
As the amalgamation property is categorical in prevarieties, these results transfer to $\mathsf{X}$ mechanically under the equivalence
in Theorem~\ref{Theorem : equivalence}.  The outcome is Theorem~\ref{Theorem : subvar with AP} below.

Here and subsequently,
$\boldsymbol{C}$ denotes the unique algebra in $\mathsf{X}$ whose lattice reduct is the chain of non-positive integers (ordered conventionally), on which
$\lozenge$ is the identity function and $\sbf=\top=0$.
If $1\leq n\in\omega$, then $\boldsymbol{C}_{2n + 1}$ is the subalgebra of $\boldsymbol{C}$ with universe $\{ - n, - n + 1, \dots, 0 \}$, while
$\boldsymbol{C}_{2n}$ is the algebra in $\mathsf{X}$ with the same Brouwerian reduct as $\sbC_{2n+1}$, but with
$\boldsymbol{f}=-1$ and $\lozenge \boldsymbol{f} = \top$ and $\lozenge a = a$ for $a \ne \boldsymbol{f}$.

\begin{Theorem}\
\label{Theorem : subvar with AP}
\begin{enumerate}
\item\label{x 1}
The nontrivial proper subvarieties of\/ $\mathsf{X}$ are\/
$\mathbb{V}(\boldsymbol{C})$\textup{,} $\mathbb{V}(\boldsymbol{C}_{n})$ for\/ \textup{$n\geq 2$,} $\mathbb{V} \{ \boldsymbol{C}, \boldsymbol{C}_{2n} \}$ for\/ $n \geq 1$ and\/ $\mathbb{V}\{ \boldsymbol{C}_{2m}, \boldsymbol{C}_{2n + 1} \}$ for\/ $n \geq m \geq 1$\textup{.}

\smallskip

\item\label{x 2}
The nontrivial subvarieties of\/ $\mathsf{X}$ with the amalgamation property are\/ $\mathsf{X}$\textup{,} $\mathbb{V}(\boldsymbol{C})$\textup{,} $\mathbb{V}(\boldsymbol{C}_{n})$ for\/ $2 \leq n \leq 4$\textup{,} $\mathbb{V} \{ \boldsymbol{C}, \boldsymbol{C}_{2n} \}$ for\/ $1 \leq n \leq 2$ and\/ $\mathbb{V}\{ \boldsymbol{C}_{2}, \boldsymbol{C}_{3} \}$\textup{.}
\end{enumerate}
\end{Theorem}

\begin{Theorem}\label{Theorem : ES for X}
Every subvariety of\/ $\mathsf{X}$ has the ES property.
\end{Theorem}

\begin{proof}
By Lemma~\ref{Lemma : X has weak ES} and Theorem~\ref{hoogland et al}, the varieties in Theorem~\ref{Theorem : subvar with AP}(\ref{x 2})
have the strong ES property.  On the other hand, finitely generated subvarieties of $\mathsf{X}$ have surjective epimorphisms, by
Theorem~\ref{expansions}(\ref{expansions 2}).  It is therefore enough, by Theorem~\ref{Theorem : subvar with AP}(\ref{x 1}),
to prove the ES property for the variety $\mathsf{K}=\mathbb{V} \{ \boldsymbol{C}, \boldsymbol{C}_{2n} \}$, where $n \geq 3$.
This is what we do now.

Whenever a member of $\mathsf{X}$ is denoted as $\sbA^+$ below, it is understood that $\sbA$ is its Brouwerian reduct.

Suppose, with a view to contradiction, that there is a $\mathsf{K}$--epimorphic inclusion $f\colon\sbB^+\mrig\sbA^+$,
with $B\neq A$.  \,As $f$ is also a $\mathsf{BRA}$--morphism, the Esakia morphism
$f_{\ast} \colon \boldsymbol{A}_{\ast} \mrig \boldsymbol{B}_{\ast}$ is not injective, by Lemma~\ref{easy}(\ref{easy 1}).
Choose distinct
$x, y \in \textup{Pr}\,\sbA$
such that $f_{\ast}(x) = f_{\ast}(y)$, i.e., $x\cap B=y\cap B$.

Let $\bup{\sbx}$
be the E-subspace of $\sbA_*$ with universe $\up{x}$.  Let $i_{x} : \bup{\sbx}
\mrig \boldsymbol{A}_{\ast}$ be
the inclusion map, so $i^{\ast}_{x}\colon{\sbA_*}^*\mrig(\bup{\sbx})^*
$ is onto, again by Lemma~\ref{easy}(\ref{easy 1}).  The canonical
isomorphism from $\sbA$ to ${\sbA_*}^*$, followed by $i_x^*$, is thus the surjective $\mathsf{BRA}$--morphism
$j_x\colon\sbA\mrig(\bup{\sbx})^*$
defined by
\begin{equation}\label{jx}
j_x(a) = \{ F \in \textup{Pr}\,\boldsymbol{A}\colon x \cup \{ a \} \subseteq F \}  \textup{ \,for all $a\in A$.}
\end{equation}
Recall that ${\boldsymbol{\mathit{Con}}}\,\sbA^+={\boldsymbol{\mathit{Con}}}\,\sbA$, so $j_x$ is
also a homomorphism from $\sbA^+$
onto an algebra $\sbE^+\in\mathsf{X}$ such that $\sbE=(\bup{\sbx})^*$.
Thus, $\sbE^+\in\mathsf{K}$.  In fact, $\sbE^+\in\mathsf{K}_\textup{FSI}$,
because ${\boldsymbol{\mathit{Con}}}\,\sbE^+={\boldsymbol{\mathit{Con}}}\,\sbE$ and $\bup{\sbx}$
has a least element (see
Lemma~\ref{easy}(\ref{Lemma : Fin Subdir Irr})).
Analogously, there exist $\sbG^+\in\mathsf{K}_\textup{FSI}$
and a surjective homomorphism
$j_{y} \colon \boldsymbol{A}^+ \mrig \boldsymbol{G}^+$, defined by
\[
j_{y}(a) = \{ F \in \textup{Pr}\, \boldsymbol{A}\colon y \cup \{ a \} \subseteq F \}  \textup{ \,for all $a\in A$,}
\]
where $\sbG$ is the dual of the E-subspace
of $\sbA_*$ on $\up{y}$.

Then $g\colon a\mapsto\langle j_x(a),j_y(a)\rangle$ is a homomorphism from $\boldsymbol{A}^+$ into $\boldsymbol{E}^+ \times \boldsymbol{G}^+$.

\begin{Claim}\label{Claim: 1 ES}
\ $g[B]\neq g[A]$.
\end{Claim}

Because $x,y\in \textup{Pr}\,\sbA$ and $x\neq y$, we may assume, by symmetry, that there exists $c \in x \smallsetminus y$. In particular,
\[
j_{x}(c) = \up{x} \text{ \ and \ } j_{y}(c) \subsetneq  \up{y}.
\]
But
$x \cap B = y \cap B$,
so for every $b\in B$,
\[
j_{x}(b) = \up{x} \text{ \ if and only if \ } j_{y}(b) =  \up{y}.
\]
Therefore, $g(c) \in g[A]\smallsetminus g[B]$, establishing Claim~\ref{Claim: 1 ES}.

Suppose $\boldsymbol{E}^+$ and $\boldsymbol{G}^+$ are both finite.  Then $g[\boldsymbol{A}^+]$ is finite.
Consider $g|_B$ as a $\mathsf{K}$--morphism from $\sbB^+$ into $g[\sbA^+]$.  Claim~\ref{Claim: 1 ES} says that this map is not onto,
so it is not a $\mathsf{K}$--epimorphism, because $\mathsf{K}$ has the weak ES property (see Lemma~\ref{Lemma : X has weak ES} and the
remark before Definition~\ref{weak es def}).
In other words, there exist $\sbD\in\mathsf{K}$ and distinct homomorphisms $h, k\colon g[\boldsymbol{A}^+]\mrig\sbD$ with $h \circ g|_B = k \circ g|_B$.
But this contradicts the fact that $f$ is a $\mathsf{K}$--epimorphism,
so we may assume, by symmetry, that $\boldsymbol{E}^+$ is infinite.

Recall that, by J\'{o}nsson's Theorem, $\mathbb{V}(\mathsf{M}\cup\mathsf{L})_\textup{FSI}=\mathsf{M}_\textup{FSI}\cup\mathsf{L}_\textup{FSI}$ for any
subvarieties $\mathsf{M}$ and $\mathsf{L}$ of a congruence distributive variety (see \cite{Jon67} or \cite{Jon95}).  In particular, as
$\mathsf{K} = \mathbb{V} \{ \boldsymbol{C}, \boldsymbol{C}_{2n} \}$,
\[
\boldsymbol{E}^+, \boldsymbol{G}^+ \in \mathbb{V}(\boldsymbol{C})_{\textup{FSI}} \cup \mathbb{V}(\boldsymbol{C}_{2n})_{\textup{FSI}}.
\]
J\'onsson's Theorem also gives $\mathbb{V}(\boldsymbol{C}_{2n})_{\textup{FSI}} \subseteq \mathbb{H}\mathbb{S}( \boldsymbol{C}_{2n} )$, so
$\mathbb{V}(\boldsymbol{C}_{2n})_{\textup{FSI}}$ consists of finite algebras, whence it excludes $\sbE^+$.  Consequently,
$\boldsymbol{E}^+ \in \mathbb{V}(\boldsymbol{C})$.

If $\boldsymbol{G}^+\in\mathbb{V}(\boldsymbol{C})$, then $g[\boldsymbol{A}^+] \in \mathbb{V}(\boldsymbol{C})$.  But then, we can apply the strong
ES property for $\mathbb{V}(\sbC)$ to Claim~\ref{Claim: 1 ES} and again contradict the fact that $f$ is a $\mathsf{K}$--epimorphism.  Therefore,
$\boldsymbol{G}^+$ belongs to $\mathbb{V}(\boldsymbol{C}_{2n})$
and is finite.  It follows that $\up{y}$ is finite.  Recall that $j_x[A]=E$.

\begin{Claim}\label{Claim : 2 ES}
\ $j_x[B]\neq E$.
\end{Claim}

As $\sbE^+$ is infinite, so is $\sbE_*$.  \,Claim~\ref{Claim : 2 ES} will be established if we can show that
$j_x[\sbB]_*$ is finite.

Observe that  $j_{x}[\boldsymbol{B}]_{\ast} =  (j_{x} f[\boldsymbol{B}])_{\ast}$. By Lemma~\ref{easy}(\ref{Lemma : Direct Images}), it is
enough to show that $(j_{x} \circ f)_{\ast}[{\boldsymbol{E}}_{\ast}]$ is finite.

The top element of $\boldsymbol{E}$ is $\up{x}$, and (\ref{jx}) shows
that ${j_{x}}^{-1} [\{ \up{x} \}] = x$, i.e., $(j_{x})_{\ast}(\{ \up{x} \}) = x$.
Because $\boldsymbol{E}$ is finitely subdirectly irreducible, $\{\up{x}\}$ is the least element of ${\boldsymbol{E}}_{\ast}$, by
Lemma~\ref{easy}(\ref{Lemma : Fin Subdir Irr}).  So, by (\ref{esakia morphism}) and the isotonicity of Esakia morphisms,
$(j_{x})_{\ast}[{\boldsymbol{E}}_{\ast}] = \bup{\sbx}$, whence the universe of
\[
(j_{x} \circ f)_{\ast}[{\boldsymbol{E}}_{\ast}] =
f_{\ast} [(j_{x})_{\ast}[{\boldsymbol{E}}_{\ast}]]
\]
is
$f_{\ast} [ \up{x} ] = \up{f_{\ast}(x)} =\up{f_*(y)}=f_{\ast} [ \up{y} ]$.
As $\up{y}$ is finite,
so is $
f_*[\up{y}]$.  This shows that
$(j_{x} \circ f)_{\ast}[{\boldsymbol{E}}_{\ast}]$ is finite, as required.

By Claim~\ref{Claim : 2 ES}, there exists $a \in A$ such that $j_{x}(a) \notin j_{x}[B]$.  As $\boldsymbol{E}^+\in\mathbb{V}(\boldsymbol{C})$, the
strong ES property for this variety implies that, for some $\sbH\in\mathbb{V}(\sbC)$, there are distinct homomorphisms
$h, k \colon \boldsymbol{E}^+ \mrig \boldsymbol{H}$  with $h|_{j_x[B]} = k|_{j_x[B]}$  and $hj_{x}(a) \neq kj_{x}(a)$.  But again, this contradicts
the assumption that $f$ is a $\mathsf{K}$--epimorphism, so $\mathsf{K}$ has the ES property.
\end{proof}

The main result of this section now follows from Theorems~\ref{Theorem : equivalence} and \ref{Theorem : ES for X}:

\begin{Theorem}\label{sugihara es}
Epimorphisms are surjective in all varieties of Sugihara monoids,
i.e., every axiomatic extension of\/ $\mathbf{RM}^\mathbf{t}$ has the infinite Beth property.
\end{Theorem}
\noindent
In contrast, even the finite Beth property fails for the weaker relevance logic
$\mathbf{R}^\mathbf{t}$ of \cite{AB75}, and for many neighbouring systems (see \cite{Urq99}
and \cite[Sec.\,4]{BH06}).

A CRL is called a \textit{positive Sugihara monoid} if it can be embedded into (the CRL-reduct of) some Sugihara monoid.
The class $\mathsf{PSM}$ of all positive Sugihara monoids is a variety \cite[Thm.\,4.2]{OR07}.  It consists of generalized Sugihara monoids in the sense of
Section~\ref{non-integral varieties section}, and it algebraizes the
negation-less fragment of $\mathbf{RM}^\mathbf{t}$.  Every
subquasivariety of $\mathsf{PSM}$ is a variety \cite[Thm.\,9.4]{OR07}.

\begin{Theorem}\label{psm}
Epimorphisms are surjective in every variety of positive Sugihara monoids,
i.e., all axiomatic extensions of the negation-less fragment of\/ $\mathbf{RM}^\mathbf{t}$ have the infinite Beth property.
\end{Theorem}

\begin{proof}
Let $\mathsf{W}$ be the class of all nuclear relative Stone algebras satisfying
\[
( \lozenge x \to x ) \lor ( ( y \lor ( y \to x ) ) \land \lozenge x ) \,\thickapprox\, \top.
\]
This variety is categorically equivalent to $\mathsf{PSM}$ \cite[Thm.\,10.4]{GR15}, and all of its subvarieties have the weak ES property \cite[Cor.\,12.5]{GR15}.
Let $\boldsymbol{D}$ and $\boldsymbol{D}_{n}$ be, respectively, the nuclear relative Stone algebra reducts of $\boldsymbol{C}$ and $\boldsymbol{C}_{n}$ for each
$n \geq 1$. From the category equivalence in \cite{GR15} and the classification of subvarieties of $\mathsf{PSM}$ that follows from \cite[Cor.\,4.3 \& 4.6]{OR07},
we can infer that the nontrivial subvarieties of $\mathsf{W}$ are
$\mathsf{W}$ itself,
$\mathbb{V}(\boldsymbol{D})$, $\mathbb{V}(\boldsymbol{D}_{n})$ for $n \geq 2$, $\mathbb{V} \{ \boldsymbol{D}, \boldsymbol{D}_{2n} \}$ for $n \geq 1$ and
$\mathbb{V}\{ \boldsymbol{D}_{2m}, \boldsymbol{D}_{2n + 1} \}$ for $n \geq m \geq 1$.  It is therefore enough to establish the ES property for these varieties.

Now $\mathsf{W}$ has the strong ES property, by \cite[Cor.\,12.5 \& Lem.\,12.8]{GR15} and Theorem~\ref{hoogland et al}.  The same is true of $\mathbb{V}(\boldsymbol{D})$,
as it is termwise equivalent to $\mathsf{RSA}$, which has the strong ES property, as noted earlier.  The finitely generated
subvarieties of $\mathsf{W}$ have the ES property, by Theorem~\ref{gsm}.
And the proof of Theorem~\ref{Theorem : ES for X} delivers (without any significant change) the ES property
for $\mathbb{V} \{ \boldsymbol{D}, \boldsymbol{D}_{2n} \}$, whenever $n \geq 1$.
\end{proof}

{\small
\noindent {\bf Acknowledgment.}  The second author thanks Miguel Campercholi for useful discussions concerning epimorphisms.  The third author thanks Zurab Janelidze for helpful comments on a partial presentation of this material.}

\end{document}